\input amssym

\magnification=\magstep1
\hsize=15.5truecm
\vsize=22truecm
\advance\voffset by 1truecm
\mathsurround=1pt

\def\chapter#1{\par\bigbreak \centerline{\bf #1}\medskip}

\def\section#1{\par\bigbreak {\bf #1}\nobreak\enspace}

\def\sqr#1#2{{\vcenter{\hrule height.#2pt
      \hbox{\vrule width.#2pt height#1pt \kern#1pt
         \vrule width.#2pt}
       \hrule height.#2pt}}}
\def\e{\epsilon}

\def\k{\kappa}
\def\o{\omega}

\def\d{\delta}

\def\n{\eta}

\def\R{\hbox{\bf R}}


\def\C{{\cal C}}



\def\th #1 #2. #3\par\par{\medbreak{\bf#1 #2.
\enspace}{\sl#3\par}\par\medbreak}
\def\rem #1 #2. #3\par{\medbreak{\bf #1 #2.
\enspace}{#3}\par\medbreak}
\def\proof{{\bf Proof}.\enspace}
\def\sqr#1#2{{\vcenter{\hrule height.#2pt
      \hbox{\vrule width.#2pt height#1pt \kern#1pt
         \vrule width.#2pt}
       \hrule height.#2pt}}}
\def\eop{\mathchoice\sqr34\sqr34\sqr{2.1}3\sqr{1.5}3}

                                                                     %
                                                                     %
\newdimen\refindent\newdimen\plusindent                              %
\newdimen\refskip\newdimen\tempindent                                %
\newdimen\extraindent                                                %
                                                                     %
                                                                     %
\def\ref#1 #2\par{\setbox0=\hbox{#1}\refindent=\wd0                  %
\plusindent=\refskip                                                 %
\extraindent=\refskip                                                %
\advance\extraindent by 30pt                                         %
\advance\plusindent by -\refindent\tempindent=\parindent 
\parindent=0pt\par\hangindent\extraindent 
{#1\hskip\plusindent #2}\parindent=\tempindent}                      %
\refskip=\parindent                                                  %
                                                                     %

\def\empty{\emptyset}

\def\raj{\restriction}

\def\nda{\mathrel{\lower0pt\hbox to 3pt{\kern3pt$\not$\hss}\downarrow}}
\def\nDa{\mathrel{\lower0pt\hbox to 3pt{\kern3pt$\not$\hss}\Downarrow}}
\def\nbot{\mathrel{\lower0pt\hbox to 4pt{\kern3pt$\not$\hss}\bot}}
\def\ekom{\mathrel{\lower3pt\hbox to 0pt{\kern3pt$\sim$\hss}\mapsto}}

\def\anR{\mathrel{\lower1pt\hbox to 2pt{\kern3pt$R$\hss}\not}}
\def\anoR{\mathrel{\lower1pt\hbox to 2pt{\kern3pt$\overline{R}$\hss}\not}}

\def\anRm{\mathrel{\lower1pt\hbox to 2pt{\kern3pt$R^{-}$\hss}\not}}

\def\ndda{\mathrel{\lower0pt\hbox to 1pt{\kern3pt$\not$\hss}\downdownarrows}}

\def\warrow{\mathrel{\lower0pt\hbox to 1pt{\kern3pt$^{w}$\hss}\rightarrow}}

\def\C{\hbox{\bf C}}

\null
\vskip 2.5truecm

\chapter{ON APPROXIMATIONS OF FEYNMAN PATH INTEGRALS}

\bigskip
\centerline{Tapani Hyttinen}

\vskip 2truecm

\chapter{Abstract}

From the point of view of mathematics, Feynman path integral approach to quantum mechanics is
very interesting because it is a mathematical enigma and at the same time a very useful tool in physics.
One (of the better) attempt(s) to build a bridge between the standard
Hilbert space operator model and
the path integral approach is based on Trotter-Kato formula. In this paper I will look how this approach looks from the point of view of
the finite dimensional approximations of the standard model
from [HH1]. In particular, I will show how to approximate Feynman path integrals and how these approximations determine the propagator
in the standard model which allows me to give a mathematically rigorous interpretation to (something like) Feynman path integral that gives correct
propagator under reasonable assumptions.

\vskip 1truecm

\chapter{1. Introduction}

In [HH1], a paper inspired by the work [Zi] of B. Zilber on model theory of quantum mechanics,
we studied finite dimensional approximations of the standard Hilbert space operator model of quantum mechanics 
and their connection to the standard model i.e. how the calculations  in the approximations determine the
results in the standard model. And if, in addition to position, also time is discrete
and the difference between two consecutive moments in time is large enough,
one has no major difficulties 
in giving Feynman path integrals a meaning in these approximations that, in the case of the free particle,
calculates correctly the propagator in the approximation itself
and which indeed looks like a (discrete) path integral
(the requirement that the difference between two consecutive moments in time is large enough comes
from divisibility questions arising
from our use of number theory in [HH1]).
In other cases this is less obvious as is the question
what is the connection  between path integrals in the approximations and in $L_{2}(\R )$, even in this
free particle case.
So it makes sense to 
look how path integrals behave in these approximations
and to study the connection between these path integrals in the approximations and
the path integrals in the continuous case. This is the topic of this paper.
The main point I will make is to demonstrate how delicate and complicated the connection between
Feynman path integrals and the discrete path integrals really is.
The use of ordinary limits does not seem to be enough, we will end up using ultralimits.
In particular the connection is much more complicated than e.g. the connection between
Riemann integrals and their obvious discrete approximations, Riemann integral is just a limit of
the discrete approximations.
The situation with Faynman path integral is also similar to
the use of Dirac's deltas in calculating the propagator: In [HH1],
we used ultraproducts to define objects that
look like Dirac's deltas in many ways, but most of them calculate the propagator incorrectly as we demonstrated in [HH1] and
in [HH2], where we continued our studies of Dirac's deltas. In fact, the situation
is even worse with path integrals,
I can not fully explain how to take the limits correctly, I just show that there is a way.
In [HH2], we were able to pinpoint when Dirac's deltas, as defined there, work.
Still I will give a mathematically rigorous interpretation to
Feynman path integral based
on these approximations and show, under rather reasonable assumptions, that the interpretation calculates
the propagator correctly.

I want to point out again that my goal here is not to give feasible ways to calculate propagators. 
Even in [HH1] this was not the goal although there we were able to calculate some propagators rather easily
using these approximations. The goal is  to try to understand  Feynman path integrals  from the point of view of mathematics.

This paper is organized as follows.
In Section 1, I will look the basics of the topic of this paper. After the basic facts about quantum mechanics,
in the first subsection, 
I will look
what Feynman path integral is.
This paper is heavily based on [HH1], and thus in the second subsection I will
recall basic notions, conventions and results from [HH1] that are needed and used in this paper.
In Section 2, I will give most of the definitions together with some remarks and, in the end, state the main theorem.
In Section 3, I will prove the theorem.

\chapter{1 Basics}

In this paper, as in [HH1],
I will study the cases in which there is only a single particle in a 1-dimensional space
and the Hamiltonian $H$, the operator for energy, is time independent. Then the time evolution operator
$K^{t}$, the solution to the Schr\"odinger equation, is $e^{-itH/\hbar}$ i.e. if we write
$\psi (t)$ for the state of the system at time $t$, then $\psi (t)=e^{-itH/\hbar}(\psi (0))$,
where $\hbar =h/2\pi$ and $h$ is the Planck constant.
Notice that now $K^{t+t'}=K^{t'}\circ K^{t}$, if $H$ depends on time, then this is not true.
However, to be really able to calculate the evolution of the system with time, one needs to know the
Feynman propagator $K(x,y,t)$ or the kernel of the time evolution operator (which should be the same if the two models
are equivalent), i.e. a function
such that for all states $\psi$,
$K^{t}(\psi )(y)=\int_{R}K(x,y,t)\psi (x)dx$ or in
Dirac notation $K(x,y,t)=<y|K^{t}|x>$ i.e. the inner product of $y>$ and $K^{t}(x>)$, where
$y>$ and $x>$ are eigenvectors of the position operator with eigenvalues $y$ and $x$, respectively
(these are called Dirac's deltas and they cannot be found from the Hilbert space in
the case of the standard model of quantum mechanics).
Since we mostly work in finite dimensional spaces, we will usually use Dirac notation (in these
finite dimensional models Dirac's deltas exist for suitable eigenvalues).
It is possible that the time evolution operator does not have a kernel/propagator but, obviously,
we will look only the cases in which it does have it.
Feynman suggested that this propagator can be calculated by using
path integrals. This method is not mathematically sound, see e.g. [Pe] or below.
[Pe] is the best introduction to Feynman path integrals that I have found and thus I will use it as my main reference to this topic 
although the paper is not published. But it seems that it has been used as
lecture notes.
The paper is easy to find from the net, from several sites, one of which is mentioned in the references below.
Also from [Me] one can find an introduction to Feynman path integrals and for further information, see the references in [Pe].
Despite of not being sound, the path integrals are an important alternative approach to quantum mechanics
(alternative to the standard model i.e. the Hilbert space operator model), 
since they have
turned out to be a very successful technique to solve problems in physics.
It follows that
the question: Can one make some mathematical sense out of Feynman path integrals?
is an interesting mathematical question. This question has been studied a lot during the last 70 years
but without much success. So I do not expect to solve it but still I think that I have something to say about
it that has not been said yet.

\chapter{1.1 Introduction to Feynman path integrals}

Let us move to the questions what Feynman path integrals are and how physicists calculate them (this will be how I see it,
some physicists may disagree). Feynman's Nobel price talk [Fe] gives historical back ground to this topic.
The most common way of describing path integral is the following: one can calculate the value of the
Faynman propagator $K(x,y,t)$ as an integral
$$A(t)\int_{P_{x,y}}e^{iS(p)/\hbar}dp,$$
where the integration is over all
paths $p\in P_{x,y}$, where $P_{x,y}$ contains the paths that satisfy the boundary condition $p(0)=x$ and $p(t)=y$,
$A(t)$ is so called normalizing constant depending on the time $t$ only
(needed to guarantee that the state has norm $1$, which also gives a way to calculate it - at least in principle)
and $S(p)$ is the classical action of the path $p$ i.e.
$$(0)\ \ \ \ \ S(p)=\int_{[0,t]}S^{*}(p(u),\dot{p}(u),u)du.$$
Here $S^{*}(p(u),\dot{p}(u),u)$ is the action density/Lagrangian of the system which in our
(non-relativistic single particle in one dimensional space) case is
$$(m\dot{p}(u)^{2}/2)-V(p(u)),$$
where $m$ is the mass of the particle and
$V$ is the function that gives the potential energy at $p(u)$.
(In general $S^{*}$ can be chosen to be anything that fits the bill, see e.g. the Lagrangian in general relativity
and for more on action,
see any book on classical mechanics or [Me].)
Notice that to calculate $S(p)$, the path needs to be differentiable (at least outside a measure zero set)
because $S(p)$ depends on the speed, $\vert\dot{p}(u)\vert$, of the particle on the path $p$ at times $u$.
However, when this method is used,
one works hard not to calculate this integral. And for an obvious reason: No one knows
what these paths are nor what is the measure used in the integration. In fact, one can prove
that the answers to these questions can not be found, see [Pe]. 

So what do the physicists do here? They split the time line from $0$ to $t$ into $n$ pieces,
each of length $\Delta t=t/n$ and
they study the expression
$$(1)\ \ \ \ \ A(t)\int_{R}...\int_{R}e^{i\Delta t(\sum_{k=1}^{n}S_{n}^{*}(x_{k-1},x_{k},\Delta t)/\hbar)}dx_{1}...dx_{n-1},$$
where $x_{0}=x$, $x_{n}=y$ and
$S_{n}^{*}(x_{k},x_{k+1},\Delta t)$ is some kind of an average Lagrangian
of the particle moving
from $x_{k}$ to $x_{k+1}$ in time $\Delta t$ i.e $\Delta tS_{n}^{*}(x_{k},x_{k+1},\Delta t)$ 
gives something like the action of the
particle moving
from $x_{k}$ to $x_{k+1}$ in time $\Delta t$. 
Still, one can not calculate this since it is not clear what $\Delta tS_{0}^{*}(x_{k-1},x_{k},\Delta t)$
should be since the action of a particle moving from $x_{k-1}$ to $x_{k}$
depends on the path the particle takes.
I have seen many attempts to solve this problem one of which is to use the classical path here
i.e. physicists replace the expression $\Delta tS_{n}^{*}(x_{k-1},x_{k},\Delta t)$ by the classical action of
the particle moving from $x_{k-1}$ to $x_{k}$ in time $\Delta t$ i.e. they look only one path here and the path is
the one classical mechanics gives. There are also many variations of this idea used in physics.
After this, the expression can be calculated, at
least in principle. Then finally, the propagator is
the limit of these when $n$ goes to infinity. Taking this limit is not easy,
I have seen infinite normalizing constants used here.

I will assume that the kernel is continuous (Assumption 4) and that the Hamiltonian is such that
the Trotter-Kato formula can be applied to the time evolution operator (Assumption 2). In the continuous case, these assumptions give
immediately (see Section 3.1) that the kernel $K(y_{0},y_{1},t)$ is
$$lim_{p\rightarrow\infty}lim_{n\rightarrow\infty}\int_{[y_{1}-1/p,y_{1}+1/p]}\int_{\R}....\int_{\R}\int_{[y_{0}-1/p,y_{0}+1/p]}$$
$$(p/2)^{2}(m/2\pi i\hbar\Delta t)^{n/2}e^{i\Delta t(\sum_{j=1}^{n}(m((x_{j-1}-x_{j})/\Delta t)^{2}-2V(x_{j-1}))/2\hbar} dx_{0}...dx_{n},$$
where $p$ and $n$ are natural numbers and $\Delta t=t/n$.
This looks a bit like (1) above. However there are three
big differences: First of all there are $lim_{p\rightarrow\infty}$, $\int_{[y_{1}-1/p,y_{1}+1/p]}$,
$\int_{[y_{0}-1/p,y_{0}+1/p]}$ and $(p/2)^{2}$ which do not exist  in (1). It seems very difficult to get rid of these,
we will take a closer look at this in Section 2 just before Assumption 4, see also the end of Section 3.1.
Then the normalizing constant $(p/2)^{2}(m/2\pi i\hbar\Delta t)^{n/2}$
or $(m/2\pi i\hbar\Delta t)^{n/2}$¨(if one considers $(p/2)^{2}$ something one should be able to get rid of)
depends not only on $t$ but also on $n$ (and $p$). Notice that the absolute value of
$(m/2\pi i\hbar\Delta t)^{n/2}$ goes to infinity when $n$ goes to infinity
(this may explain the infinite normalizing constant) and so one can not
pull it out of the limit. Also the term itself oscillates
wildly.
Finally, $\Delta t(m((x_{j-1}-x_{j})/\Delta t)^{2}-2V(x_{j-1}))/2$ may be much smaller than
the action of any path the particle can take moving from  $x_{j-1}$ to $x_{j}$ in time $\Delta t$.
Still my starting point will be this observation made in the approximations when I try to understand
(0) above i.e. the first version Feynman path integrals, which is the version that I am interested in.

\chapter{1.2 Introduction to [HH1]}

In this subsection I describe the basic constructions and results from [HH1]. Later I will recall more results from [HH1]
when they are needed. Eventually we will end up recalling all the results from [HH1] that one needs to know
to read this paper.

I start by explaining what the approximations 
$H_{N}$, $N$ a natural number, of the standard model from [HH1] are: 
$H_{N}$ is a Hilbert space with an orthonormal basis
$u(n)$, $n<N$, and each $u(n)$ is an eigenvector of the position operator $Q_{N}$ with eigenvalue $n/\sqrt{N}$ if
$n<N/2$ and otherwise the eigenvalue is $(n-N)/\sqrt{N}$. We are interested in only those $N$ s.t. $\sqrt{N}$ is a 
natural number and divisible by all small natural numbers i.e. we require that the
set $\{ N<\o\vert\ \sqrt{N}\ \hbox{\rm is a natural number divisible by}\ n\}$ belongs to the ultrafilter $D$ for all
natural numbers $n$ and here $D$ is the filter that was used in [HH1] to build a connection
between the standard model and our approximations, see below.
As usually, if something happens in a set that belongs to $D$, we say that it happens for almost all $N$.
Then as in [HH1], we think that a claim is true if it is true 
for almost all $N$ or
in fact true for almost all $N_{k}$ or for almost all $N^{d}_{k}$, see below
(we will work with several ultrafilters $D$) and with definitions we do the same
i.e. it is enough for the definitions to make sense only for almost all $N$.
We write $X_{N}$ for the set of $Q_{N}$-eigenvalues in $H_{N}$ and for $x\in X_{N}$, we denote by $x>$,
following Dirac notation, the eigenvector $u(n)$
with the eigenvalue $x$. The momentum operator $P_{N}$ is got essentially by conjugating $Q_{N}$ with the (discrete) Fourier transform i.e.
$P_{N}=h\Phi_{N}\circ Q_{N}\circ \Phi_{N}^{-1}$ where $\Phi_{N}$, the Fourier transform, is a unitary operator (isometric automorphism)
from $H_{N}$ to $H_{N}$ and $h$ is the Plank constant as before.
In particular, letting $q=q_{N}=e^{i2\pi/N}$ the vectors
$$v(n)=\Phi_{N}(u(n))=N^{-1/2}\sum_{m=0}^{N-1}q^{nm}u(m),$$
are $P_{N}$-eigenvectors and they also form a basis for $H_{N}$. Notice that this determines what
the Fourier transform $\Phi_{N}$ is.
Then we picked a suitable ultrafilter $D$ (see above) on natural numbers and took a metric ultraproduct of the approximations.
More precisely, we took the ultraproduct of the approximations, got rid of all elements with infinite norm and finally
equalized vectors according to the equivalence relation of being infinitely close to each other. In [HH2] we looked also other versions of this construction.
In this ultraproduct, the operators $Q_{N}$ and $P_{N}$ do not have a uniform modulus of continuity
(they approximate unbounded operators) and thus their (metric) ultraproducts are not
well-defined in the metric ultraproduct. However, they are well defined (as unbounded self adjoint operators, for details see [HH1])
in a subspace, call it $H_{u}$ and we call the
operators in $H_{u}$ simply $Q$ and $P$,
and, in addition, we can choose $H_{u}$ so that there is an isometric isomorphism $F:L_{2}(\R )\rightarrow H_{u}$
that maps the usual position and momentum operators (also called $Q$ and $P$) to $Q$ and $P$. In addition, we showed how calculations in the
spaces $H_{N}$ can be transferred to results in $H_{u}$ and thus results in $L_{2}(\R )$ (this is more complicated than what
one might think).

However, in this paper we look operators like $e^{-itP^{2}/2m\hbar}$ which are unitary and thus bounded and so
the problems we had in [HH1] with $P$ and $Q$ do not arise here, we have already dealt all these problems in [HH1]. Also to a large extend,
in [HH1] we already figured out what operators the ultraproducts give us,
but see Question 1, Assumption 1
and the paragraph just before Assumption 1 below.

For later use, let us look more closely how these $F$ and $H_{u}$ were found. For e.g. compactly supported smooth functions
$f\in L_{2}(\R )$ (as usually we overlook the fact that the elements of $L_{2}(\R )$ are actually equivalence classes)
we define
$$(2)\ \ \ \ \ F_{N}(f)= N^{-1/4}\sum_{x\in X_{N}}f(x)\vert x>\in H_{N},$$
where $f(x)\vert x>$ means the vector one gets by multiplying $x>$ by $f(x)$.
Then for these $f$ we let
$$(3)\ \ \ \ \ F(f)=(F_{N}(f))_{N<\o}/D/\sim ,$$
where $\sim$ is the equivalence relation being infinitely close to each other. This determines a unique
isometric function $F$ from $L_{2}(\R )$ to the metric ultraproduct and we let $H_{u}$ to be the image of this $F$.
Then many things from $L_{2}(\R )$ can be calculated using these approximations $H_{N}$.
Suppose that $U$ is a unitary operator on $L_{2}(\R )$
and $U_{N}$ are unitary operators in the spaces $H_{N}$ such that they approximate $U$ i.e. their metric ultraproduct restricted to
$H_{u}$ is the image of $U$ under $F$ (in particular, the ultraproduct maps $H_{u}$ to $H_{u}$).
Then clearly, by (3) and \L os's theorem,
for compactly supported
smooth functions $f,g\in L_{2}(\R )$,
$$(4)\ \ \ \ \ <g\vert U\vert f>=lim_{D}<F_{N}(g)\vert U_{N}\vert F_{N}(f)>,$$
where 
$lim_{D}h(N)$ for $h:\o\rightarrow\C$, is the unique complex number $q$ such that for all
$\e >0$,
$$\{ N<\o\vert\ \vert h(N)-q\vert <\e\}\in D,$$ if such $q$ exists (by compactness, it exists if the values of $h$
are bounded which is the case every time we use this notation). Finally, since $U$ was assumed
to be unitary (convergences are preserved) it is easy to see that (3) and (4) hold also for
all functions $f\in L_{2}(\R )$ that can be approximated in a nice regular way by compactly supported smooth functions
when $F_{N}$ for these is defined as in (2),
in particularly this is true for step functions from $L_{2}(\R )$ with finitely many steps (these are compactly supported).
More precisely, if $f:\R\rightarrow\C$ is such that there are compactly supported smooth functions $f_{i}$, $i<\o$,
such that the sequence converges to $f$ and $limsup_{N<\o}\vert F_{N}(f)-F_{N}(f_{i})\vert_{2}$ goes to zero when $i$ goes to infinity,
then $F(f)=(F_{N}(f))_{N<\o}/D/\sim$ since $F(f)$ is the limit of ($F(f_{i}))_{i<\o}$.
Also (4) holds  since (letting $(g_{i})_{i<\o}$ be for $g$ as $(f_{i})_{i<\o}$ was for $f$)
$U(f)$ is the limit of $(U(f_{i}))_{i<\o}$
and $<F_{N}(g_{i})\vert U_{N}\vert F_{N}(f_{i})>$ approximate $<F_{N}(g)\vert U_{N}\vert  F_{N}(f)>$ as closely as we want when $N$
and $i$ are large enough ($U_{N}$ preserves the distance $\vert F_{N}(f)-F_{N}(f_{i})\vert_{2}$).
This was shown and used in [HH1] in the case of time evolution operators.
Notice that this is not true for all $f\in L_{2}(\R )$: Let $f:\R\rightarrow\C$ be such that $f(x)=0$ if
$x$ is a rational number or  $\vert x\vert >1$ and otherwise $f(x)=1$. Then
$f\in L_{2}(\R )$ but $F(f)\ne (F_{N}(f))_{N<\o}/D/\sim$ if $F_{N}(f)$ is defined as in (2) above.

Finally, in [HH1] we used number theory to prove and calculate our results. It followed that we need to
scale the units so that some values are natural number and some rational numbers. In this paper we use the scaling
for the free particle from [HH1] i.e. scaling we needed in the case $H=P^{2}/2m$. In [HH1], in the case of the harmonic oscillator,
a different scaling was needed. However, in this paper we deal with these more complicated Hamiltonians by using
Trotter-Kato formula and thus this scaling for the free particle works everywhere.

We fix the time $t>0$ and the mass $m$ of the particle.
Then we need
to scale the units so 
that $th/2m$ is an integer. We can get this by scaling e.g. only the length. However, later when  we start splitting the time we are in problems,
we need, at least eventually, to rescale the units and end up studying questions in a smaller scale (this may have a connection  to what
physicists study in QFT
under the name renomalization group, brave readers can try to learn about this e.g from [Me]). We will look this later. 

\medskip

{\bf Remark 1}. Scaling of the length does not affect to the position operator (only to our interpretation of the numbers). However,
this is not true with the momentum operator, it depends on the numerical value of $h$ and this depends on the scaling.

\medskip

We need also that $t/\pi$ and $th/\pi$ are rational numbers. 
We can scale units for time and mass so that $t/\pi$ and $t/m$ are a rational numbers and then the scaling
of length must make also $h$ a rational number (the unit of $h$
is $mass\times length^{2}/time$). After this it is easy to see that we can
choose the scaling for length so that, in addition, $th/2m$ is any positive natural number we want.
To simplify not only our formulas but also our constructions,
we choose these scalings so that $th/2m=1$ i.e. $th/m=2$.
Then $2$ will appear a lot in our formulas and most of the time it is this $th/m$.
To some extend we will keep track of this $th/m$. We will return to these scalings, in the next section.

Finally, we aim to calculate $K(y_{0},y_{1},t)$ in Feynman path integral style and here our techniques
requires also that $y_{0}$ and $y_{1}$ are rational numbers (in the chosen unit of length but notice that the rational numbers are dense in $\R$
and we will assume that for our fixed $t$,
$K(y_{0},y_{1},t)$ is a continuous function $\R^{2}\rightarrow\C$ and thus this is not really a restriction).

\chapter{2 Definitions and the main theorem}

We start by looking the case of the free particle i.e. the case when $H=P^{2}/2m$. Here the Lagrangian and the classical Hamiltonian
are the same and thus this is not the most interesting case in it's own
(although to demonstrate the points I am making, this would be have been good enough for me)
but we will need it when we look more general cases.

We start from the question: what the approximations 
of the time evolution operator are. In this free particle case this is easy. In general it is not (we will return to this later):
E.g. in [HH1], in the case of the harmonic oscillator we used rather tricky approach for two reasons.
First, it is not clear that the natural choices actually approximate the time evolution operator in $L_{2}(\R )$
and second, even if they do,
the propagators for these natural choices
would not be easy to calculate in the approximations. But in this free particle case
the time evolution operator $K^{t}_{N}$ is $e^{-itP_{N}^{2}/2m\hbar}$ in the approximation $H_{N}$. 
In [HH1] we showed that the ultraproduct of these operators gives the right operators in
the image $H_{u}$ of the embedding $F$ of $L_{2}(\R )$ to the ultraproduct. 
In the approximation $H_{N}$, I will write simply $K^{t}$ for the approximation $K^{t}_{N}$ and similarly $P$ and $Q$ for $P_{N}$ and $Q_{N}$,
when $N$ is clear from the context (this makes formulas more readable).

For all $k<\o$ and functions $\d:\o\rightarrow\o$ (we will choose this $\d$ later - in fact many times) we let $N_{k}=\Pi_{i=0}^{k}p_{i}^{2k}$
and $N_{k}^{d}=2^{2\d (k)}N_{k}$, where $p_{k}$ is the $(k+1)$st prime number ($p_{0}=2$ etc.).
We write $H_{k}=H_{N_{k}}$ (there will be no confusion for the two possible meanings of $H_{k}$).
$H^{d}_{k}$ is defined a bit differently. The idea is that in $H^{d}_{k}$ we have changed the scaling of length
so that one unit of length in $H_{k}$ is $2^{\d (k)}$ units of length in $H^{d}_{k}$. This means that
we have build $H^{d}_{k}$ exactly as $H_{N^{d}_{k}}$ was build only as $h$ we have used $2^{2\d (k)}h$ (we recall that the unit of
$h$ is $mass\times length^{2}/time$). Notice that with this $h$, e.g. $th/\pi$ is still a rational number, see Section 1.2, and
$\sqrt{N_{k}^{d}}$ and $\sqrt{N_{k}}$ are integers (divisible by many natural number).
We write $X_{k}$ for $X_{N_{k}}$ and
$X^{d}_{k}$ for the set of $Q_{N^{d}_{k}}$-eigenvalues and keep in mind
that the units in $H^{d}_{k}$ are not the same as those in $H_{k}$.

Lemma 1 below (from [HH1]), will be in an important role in this paper. It also gives an interesting connection between $H_{k}$ and $H^{d}_{k}$
which we will look first.
We define an embedding $G_{k}:H_{k}\rightarrow H^{d}_{k}$ so that $u(n)\in H_{k}$ is mapped to $u(2^{2\d (k)}n)\in H^{d}_{k}$
i.e. $G_{k}$ maps $x>$ to $(2^{\d (k)}x)>$ and notice that $2^{\d (k)}x$ in units of $H^{d}_{k}$ is $x$ in units of $H_{k}$.
Thus if we switch to use units of $H_{k}$ also in $H^{d}_{k}$, $G$ maps $x>$ to $x>$ (in the proof of the main theorem,
we will  switch from units of $H^{d}_{k}$ back to the units of $H_{\k}$ in some formulas originally calculated
in $H^{d}_{\k}$ using the units of $H^{d}_{k}$).
The function $G_{k}$ is clearly an isometric embedding and the following lemma from [HH1] tells what it does to the time evolution operator
in the case of the free particle. Also
there is a reason why the formula for the propagator is written in such an obscure way, $K^{*}(x_{0},x_{1},t)$ is the correct value
of the propagator in $L_{2}(\R )$ in this free particle case. Before giving the lemma we want to point out that embeddings like
$G_{k}$ are not the natural way of embedding approximations to better approximations e.g. if one takes the direct/ultra limit of all approximations under 
embeddings of this kind, the result is not anything like the standard model. But these embeddings, surprisingly, will give a nice connection
between the two approximations.

\th 1 Lemma. ([HH1]) In the free particle case in $H_{N}$, if $th/m$ divides $\sqrt{N}(x_{1}-x_{0})$,
$x_{0},x_{1}\in X_{N}$, then
$<x_{1}\vert K^{t}\vert x_{0}>=N^{-1/2}thm^{-1}K^{*}(x_{0},x_{1},t)$ and otherwise it is $0$, where
$K^{*}(x_{0},x_{i},t)=(m/2\pi i\hbar t)^{1/2}e^{im(x_{1}-x_{0})^{2}/2\hbar t}$.

So the range of $G_{k}$ is closed under $K_{N^{d}_{k}}^{t}$ and $G_{k}\circ K^{t}_{N_{k}}\circ G_{k}^{-1}=K^{t}_{N^{d}_{k}}\raj rng(G_{k})$
since $\sqrt{2^{2\d (k)}h}/\sqrt{N^{d}_{k}}=\sqrt{h}/\sqrt{N_{k}}$, the formula in the exponent is dimensionless and
for $x_{0},x_{1}\in X_{k}$,
$t2^{2\d (k)}h/m$ divides $\sqrt{N^{d}_{k}}(2^{\d (k)}x_{1}-2^{\d (k)}x_{0})$ iff $th/m$ divides $\sqrt{N_{k}}(x_{1}-x_{0})$
and for $x_{0},x_{1}\in X^{d}_{k}$,
if one of $x_{0}>$ and $x_{1}>$ is in $rng(G_{k})$ and the other one is not, then $t2^{2\d (k)}h/m$
does not divide $\sqrt{N^{d}_{k}}(x_{1}-x_{0})$ and so the propagator will be $0$.
However, notice that the propagator has the unit $1/length$ in $L_{2}(\R )$ (units are in the approximations even trickier than
in $L_{2}(\R )$, we will take a closer look at this in Section 3.1) and  so although the  numerical value
of
$$<2^{\d (k)}x_{1}\vert K^{t}_{H^{d}_{k}}\vert 2^{\d (k)}x_{0}>$$
in $H^{d}_{k}$ is the same as the numerical value of
$<x_{1}\vert K^{t}_{H_{k}}\vert x_{0}>$ in $H_{k}$, with the unit the value in
$H^{d}_{k}$ is $2^{\d (k)}$ times bigger than in $H_{k}$ (see above).
This is as it should, basically because of the spiking
phenomenon (i.e. 'otherwise it is $0$'-part in
Lemma 1), see [HH1]. In Subsection 3.1 we will take a closer look at effects of our changes in the scalings. 

We would like to have that for $n=2^{2\d (k)}$ the range of $G_{k}$ is closed under $K_{N^{d}_{k}}^{t/n}$
and $G_{k}\circ K^{t/n}_{N_{k}}\circ G_{k}^{-1}=K^{t/n}_{N^{d}_{k}}\raj rng(G_{k})$. 
However, this is not true although I believe that the operators approximate each other in
a reasonable and useful way.

In [HH1] we showed that the approximations of the operators studied in the paper indeed approximate the operators.
However, we were not able to show that always the obvious approximations approximate the operators correctly
(i.e. as in Lemma 2 below).
As mentioned above, we were not able to show even in the harmonic oscillator case that the operators
$e^{-it((P^{2}_{N}/2m)+V_{N})/\hbar}$  approximate $e^{-it((P^{2}/2m)+V)/\hbar}$ from  $L_{2}(\R )$ correctly
(for potential energy operators
$V_{N}$ and $V$, see below
and in this harmonic oscillator case we were able to find the approximations but they were not these
obvious ones).
In the (lucky?) cases when the operators $e^{-it((P^{2}_{N}/2m)+V_{N})/\hbar}$ actually  approximate $e^{-it((P^{2}/2m)+V)/\hbar}$ correctly,
we might be able to get much smoother theory than the one
below because
e.g. we would not need Trotter formula in $L_{2}(\R )$ (i.e. Assumption 2 below), we just apply it in the spaces $H_{N}$ and there it holds always
since the spaces have finite dimensions. In fact in the spaces $H_{N}$ the version of Trotter formula we have is much stronger than
the one we will use in $L_{2}(\R )$.
However, there may also be problems.  The needed growth pace of $\d$ might be too much for Remark 2 below and also
the problem with $K^{t/n}_{N^{d}_{k}}$ mentioned in the previous paragraph may couse problems if one tries
to get rid of the two step approach we use in Section 3,
where we look first constant functions $\d$ and then by diagonalizing find the final $\d$.

However, the following approximation result was proved in [HH1] (Proposition 2.15):

\th Lemma 2. For all natural numbers $n>0$,
the operators $e^{-itP^{2}_{N}/2m\hbar n}$ approximate the operator $e^{-itP^{2}/2m\hbar n}$ from  $L_{2}(\R )$ i.e.
for all ultrafilters $D$ as in Section 1.2, the metric ultraproduct of operators $e^{-itP^{2}_{N}/2m\hbar n}$
restricted to $H_{u}=rng(F)$ is the same as the image of $e^{-itP^{2}/2m\hbar n}$ in the embedding $F$.

Now to the potential energy operator $V$.
We pick a function $f_{v}:\R\rightarrow\R$ so that $V(\phi )(x)=f_{v}(x)\phi (x)$ is a self adjoint operator
in $L_{2}(\R )$ (in particular, it is densely defined e.g. $f_{v}(x)=kx^{2}/2$ where $k$ is a positive real, a
string constant). We approximated $V$ in $H_{N}$ the obvious way: $V_{N}(u(n))=f_{v}(x)u(n)$ where
$x$ is the $Q$-eigenvalue of $u(n)$ and if it is clear from the context we also write $V$ for $V_{N}$.
In $H_{k}^{d}$ we need to take into account the scaling: $V_{N_{k}^{d}}(u(n))=2^{2\d (k)}f_{v}(x/2^{\d (k)})u(n)=f_{v}^{d}(x)u(v)$, where
$x$ is the $Q_{N_{k}^{d}}$-eigenvalue of $u(n)$.
We write $V$ also for $V_{N^{d}_{k}}$ when there is no risk of confusion.
And again keep in mind that the unit of $f_{v}(x)$ is that of energy (it should give the potential energy).

In the cases when $H=(P^{2}/2m)+V$, we use Trotter-Kato formula to calculate the propagator i.e.
we will, in a sense, split $H$ into two operators (Hamiltonians) $P^{2}/2m$ and $V$. So we need to know
what happens in the case $H=V$. We start by looking the embeddings $G_{k}$.
For all $n>0$, if $H=V$, then
the range of $G_{k}$ is obviously closed under $K^{t/n}_{N_{k}^{d}}$ and 
$G_{k}$ maps $K^{t/n}$ in $H_{N_{k}}$ to $K^{t/n}\raj rng(G_{k})$ in $H_{k}^{d}$ (in the definition of $K^{t}$, $V$ is divided by $\hbar$,
which in the units of $H_{k}^{d}$ is $2^{2\d (k)}$ times it's value in the units of $H_{k}$).
Also in [HH1], we pointed out (Proposition 2.16) that
for all natural numbers $n>0$,
the operators $e^{itV_{N}/\hbar n}$ approximate the operator $e^{itV/\hbar n}$ from $L_{2}(\R )$
as in Lemma 2 in the case of harmonic oscillator i.e. when $V$ is as in the harmonic oscillator.
In fact this is much easier to prove than
Lemma 2
and the proof generalizes to many functions $f_{v}$. But not for all and thus we need to assume:

\medskip

{\bf Assumption 1}: The function $f_{v}$ is such that for all natural numbers $n>0$,
the operators $e^{-itV_{N}/\hbar  n}$ approximate $e^{-itV/\hbar n}$ from  $L_{2}(\R )$ as in Lemma 2.

\medskip

In fact we use Assumption 1 to conclude,
in the case when $\d (k)=s$ for all $k<\o$,
that the operators $e^{-itV_{N^{d}_{k}}/\hbar n}$ approximate the operator
$e^{-itV^{s}/\hbar n}$ from $L^{s}_{2}(\R )$, see Section 3.1. The difference between these two cases is in the choice of
the units and this property does not depend on this choice. The same is true in the case of Lemma 2.

Now we are interested in calculating the propagator for the time evolution operator in the case when $H=(P^{2}/2m)+V$.
We want to do this in a path integral style. As pointed out above for this we want to use Trotter formula
(Trotter-Kato formula, Lie product formula). 
For this we need to make the following assumption:

\medskip

{\bf Assumption 2}: In $L_{2}(\R )$, the operator $e^{-it((P^{2}/2m)+V)/\hbar}$ is the strong limit of the operators
$(e^{-itP^{2}/2m\hbar n}e^{-itV/\hbar n})^{n}$, when $n$ goes to infinity (strong limit = pointwise convergence).

\medskip

Theorem VIII.30 from [RS] gives the classical result on when Assumption 2 holds and from it, it follows that
Assumption 2 holds for the harmonic oscillator (it holds trivially for the free particle).
But if I have understood it correctly, e.g. if $f_{v}$ gives the potential
energy of Coulomb force, Assumption 2 fails (it has been challenging to find from the literature what is known
about Assumption 2, I have asked people that work in functional analysis, and physicists use it routinely).
However, here we are not interested in these question, we are
trying to understand Feynman path integrals at least in some cases no matter how simple.

{\bf Question 1}: Do Assumptions 1 and 2 imply that the operators
$$e^{-it((P^{2}_{N}/2m)+V_{N})/\hbar}$$
approximate the operator $e^{-it((P^{2}/2m)+V)/\hbar}$ from  $L_{2}(\R )$ correctly (i.e. as in Lemma 2 and
at least for some ultrafilter $D$)?
Keeping in mind that Trotter formula holds (strongly) in the approximations,
it looks at first that what we do below implies that the answer is yes. However, there is a problem.
We do not seem to have enough control on our final function $\d =\d^{*}$ (see below) to push the proof through.
However, I believe that a closer inspection of the arguments below give a positive answer to this question.

For $r\le 2^{2\d (k)}$ and $\Delta t=\Delta t(\d ,k)=t/2^{2\d (k)}$ we write
$$L_{k}^{r\Delta t}=(e^{-itP^{2}_{N^{d}_{k}}/2^{2\d (k)}2m\hbar}e^{-itV_{N^{d}_{k}}/2^{2\d (k)}\hbar})^{r}$$
and thus
$$L_{k}^{t}=(e^{-itP^{2}_{N^{d}_{k}}/2^{2\d (k)}2m\hbar}e^{-itV_{N^{d}_{k}}/2^{2\d (k)}\hbar})^{2^{2\d (k)}}.$$
So this $L_{k}^{t}$ depends on $\d (k)$. In the last subsection, Subsection 3.2, it is not always clear what $\d (k)$ is.
Thus we write also $L^{t}_{ks}$ and by this we mean $L^{t}_{k}$ is the case when $\d (k)=s$.
Then Lemma 2 and Assumption 1 give immediately:

{\bf Remark 2}. If $\d$ is a constant function, say $\d (n)=s$, then the operators $L_{k}^{t}$
approximate the operator
$$(e^{-itP^{2}/2^{2s}2m\hbar}e^{-itV/2^{2s}\hbar})^{2^{2s}}$$
from $L_{2}(\R )$ as in Lemma 2 assuming that for all $p<\o$, $\{ N^{d}_{k}\vert\ k>p\}\in D$
and that the standard model $L_{2}(\R )$ is build using units that are used in the spaces $H_{N^{d}_{k}}$
(in particular, the numerical value of $h$ is the same in all spaces,
and notice also that if for all $p<\o$, $\{ N^{d}_{k}\vert\ k>p\}\in D$, then for all
$p<\o$
$$\{ N<\o\vert\ \sqrt{N}\ \hbox{\sl is a natural number and divisible by}\ p\}\in D).$$
Later we will write
$L^{s}_{2}(\R )$ for this standard model build using the units that are used in the spaces $H_{N^{d}_{k}}$
when $\d (k)=s$.

We write $n^{*}=n^{*}(k)=n^{*}(\d ,k)=2^{2\d (k)}$ and
$$S^{*}(x,y,t)=(m(y-x)^{2}/2t^{2})-f_{v}^{d}(x).$$
When $S^{*}$ is used, it is calculated in $H^{d}_{k}$ and thus, a priori, one uses the units from $H^{d}_{k}$.
However, each time we use it in formulas for propagators, it is part of an expression that as whole is dimensionless and we will end up using
the original units (from $H_{k}$). We will also end up using the original units for the elements of $X^{d}_{k}$
and for the same reason and, of course, because writing formulas using different units is most inconvenient.
We also write $X^{2d}_{k}$ for the set of all $x\in X^{d}_{k}$ such that $x/2\in X^{d}_{k}$
i.e. $X^{2d}_{k}$ is the set of all $x\in X^{d}_{k}$ such that $2$ divides $\sqrt{N^{d}_{k}}\ x$, cf. Lemma 1 above
(so again this $2$ is $th/m$ - as calculated in the original units). 

We will take sums that run over all elements of $X^{*d}_{k}=(X^{2d}_{k})^{n^{*}-1}\cup (X^{d}_{k}-X^{2d}_{k})^{n^{*}-1}$. One can get nicer formula if
one uses integral notation here, one can bury complicated normalizing constants into the measure.
Also our aim is to understand something written using integrals which typically can be seen as limits
of discrete sums and so it is natural to use integral notation also in this discrete case.
Thus we define a measure $\mu =\mu^{k}$ to $X^{*d}_{k}$ so that it is the counting measure normalized so that
$dx_{1}...x_{n^{*}-1}=\sqrt{2}(N^{d}_{k}/2)^{-n^{*}/2}$
(here $2$ is again $th/m$ in two of the three places)
i.e.
$$\mu (X^{*d}_{k})=2\vert X^{2d}_{k}\vert^{n^{*}-1}\sqrt{2}(N^{d}_{k}/2)^{-n^{*}/2}.$$

If one looks our main theorems, Theorem 1 and Lemma 5 below, one notices that there the integral
is multiplied by $2^{-2\d (k)}N_{k}^{-1/2}$ (the measure we use there is
based on this measure). So we probably should multiply our measure also by this 
(we have not done this for technical reasons). 
However non of this affects
to our next consideration (or it does but it makes the convergence stronger).

Suppose that $Y_{k}$ is the set of all sequences $(x_{1},...,s_{n^{*}-1})\in X^{*d}_{k}$
such that $0\le x_{i}<1$ for all $0<i<n^{*}$. Then
$$lim_{k\rightarrow\infty}\mu^{k}(Y_{k})=lim_{k\rightarrow\infty}2(\sqrt{N^{d}_{k}}/2)^{(n^{*}-1)}\sqrt{2}(N^{d}_{k}/2)^{-n^{*}/2}=0.$$
The same happens for all integers $n<m$ in place of $0<1$. Since below we will show that our measure works
in our path integral, this may indicate the nature of the problem of finding
a Lebesgue/Riemann style path integral in the continuous case (i.e. in $L_{2}(\R )$) and why infinite normalizing constants are some times needed
when one tries to calculate Feynman path integral in the continuous case.
One expects that e.g. if we let $P_{00}^{nm}$, $n<0<m$, be the set of all continuously differentiable paths $p$
for which $n\le p(u) <m$ for all $0<u<t$, $p(0)=0$ and
$p(t)=0$, then $P_{00}^{nm}$ is measurable. But now the calculation above gives measure zero to the
set of approximations of the paths of this set at the limit i.e. in the continuous case. For the approximations of
these paths, see below.

Before giving further definitions let me make an observation that concerns the expression
$$\Delta t(\sum_{j=1}^{n^{*}}S^{*}(x_{j-1},x_{j},\Delta t)),$$
keep in mind that $\Delta t=t/2^{2\d (k)}$
This will motivate the rest of the definitions in this section.

\medskip

{\bf Remark 3}. Suppose that we use in the sets $X^{2d}_{k}$ the same units as in the sets $X_{k}$
i.e. they become sets $\{ 2^{-\d (k)}x\vert\ x\in X^{2d}_{k}\}$ and that we also write the values of
$f^{d}_{v}$ in units of $H_{k}$ taking in the account that we use units of $H_{k}$ in $X^{2d}_{k}$ i.e. $f^{d}_{v}$
becomes $f_{v}$ (in the definition of
$S^{*}$).
Suppose also that $p:[0,t]\rightarrow\R$
is a continuously differentiable path (i.e. velocity is continuous). Now for all $k<\o$ and $0\le j\le n^{*}$,
let $p^{k}_{j}$ be the largest element in $X^{2d}_{k}$ that is $\le p(jt/n^{*})$
($min X^{2d}_{k}$ if there is no such element). 
Now one notices,
assuming that $f_{v}$ is continuous, that the classical action
of the path $p$ i.e.
$$\int_{[0,t]}((m\dot{p}(u)^{2}/2)-f_{v}(p(u)))du,$$
is
$$lim_{k\rightarrow\infty}\Delta t(\sum_{j=1}^{n^{*}(k)}S^{*}(p_{j-1},p_{j},\Delta t)),$$
if $\d$ is non-decreasing, $lim_{k\rightarrow\infty}\d (k)=\infty$
and $2^{2\d (k)}/\sqrt{N^{d}_{k}}$ goes to zero when $k$ goes to infinity. This last assumption holds
for all $\d$'s that we will look, in fact, for them $\d (k)\le k$ for all $k$ large enough.
Before proving this, let us look some variations of this result.

It follows that we can weaken the first two assumption on  $\d$, if we switch to $lim_{D^{*}}$
for an ultrafilter $D^{*}$ on $\o$ i.e.
if for all $q <\o$,
$\{ k<\o\vert\ \d (k)>q\}\in D^{*}$, then the action is
$$lim_{D^{*}}\Delta t(\sum_{j=1}^{n^{*}(k)}S^{*}(p_{j-1},p_{j},\Delta t)).$$

Also we can take the limit in two steps (and get rid of the annoying assumption
that $2^{2\d (k)}/\sqrt{N_{k}}$ goes to zero when $k$ goes to infinity):
For non-zero $s<\o$,
let $\d =\d_{s}$ be the constant function $\d_{s}(n)=s$ for all $n<\o$
and write
$$S_{s}(p)=lim_{k\rightarrow\infty}\Delta t(\sum_{j=1}^{n^{*}(\d_{s},k)}S^{*}(p_{j-1},p_{j},\Delta t)).$$
Then the action is $lim_{s\rightarrow\infty}S_{s}(p)$.

So one can think the sequence $((p^{k}_{j})_{j\le n^{*}(k)})_{k<\o}$ as an approximation of
the path $p$ in the strong sense that it does not only approximate the path itself in the geometric sense
but it approximates also the action of the path.
However, the problem here is that this (or variations of this) do not live in the ultraproduct.
We will return to this later.

\medskip

{\bf Proof of Remark 3}. Since $f_{v}$, $p$ and $\dot{p}$ are continuous,
Lagrangian is 
Riemann integrable and since 
$$(p((j-1)t/n^{*})-p(jt/n^{*}))^{2}/\Delta t^{2}$$ is
a value $\dot{p}^{2}$ gets at some time in $[(j-1)t/n^{*},jt/n^{*}[$ (by the mean value theorem),
the action is the limit of
$$\Delta t\sum_{j=1}^{n^{*}}(m(p((j-1)t/n^{*})-p(jt/n^{*}))^{2}/2\Delta t^{2})-f_{v}(p((j-1)t/n^{*}))$$
when $k$ goes to infinity.
So it is enough to show that the difference between this expression and the expression from the claim goes to $0$, when $k$ goes to infinity.

We start by noticing that
$\vert p(jt/n^{*})-p^{k}_{j}\vert\le 2(N^{d}_{k})^{-1/2}$ for large enough $k$ (the range of $p$ is bounded).
Now let $c$ be an upper bound for the absolute values of the values of $\dot{p}$.
Then
$$\vert p((j-1)t/n^{*})-p(jt/n^{*})\vert\le c\Delta t.$$
Also given $\epsilon >0$ for all $k$
large enough, $\vert f_{v}(p((j-1)t/n^{*}))-f_{v}(p^{k}_{j-1})\vert <\epsilon$ since $f_{v}$ is uniformly continuous
in a closed interval that contains $rng(p)$.
With these, for large enough $k$, first
$$\vert (p((j-1)t/n^{*})-p(jt/n^{*}))^{2}-(p^{k}_{j-1}-p^{k}_{j})^{2}\vert\le $$
$$(4\vert p((j-1)t/n^{*})-p(jt/n^{*})\vert /\sqrt{N^{d}_{k}})+4(N^{d}_{k})^{-1}\le (4c\Delta t/\sqrt{N^{d}_{k}})+4(N^{d}_{k})^{-1}.$$
Thus the difference between the two expressions is at most
$$(n^{*}(k)m2c/\sqrt{N^{d}_{k}})+2n^{*}(k)m(\Delta tN^{d}_{k})^{-1}+n^{*}(k)\Delta t\epsilon .$$
The first of these terms go to zero because of our assumption that $2^{2\d(k)}/\sqrt{N^{d}_{k}}$ goes to zero.
Similarly the second term goes to zero since if $2^{2\d(k)}/\sqrt{N^{d}_{k}}$ goes to zero, also
$(2^{2\d(k )})^{2}/N^{d}_{k}$ goes to zero.
The last term goes to
$t\epsilon$ (actually it is $t\epsilon$). So the difference is at most $t\epsilon$ and since this holds for all
$\epsilon >0$, the difference, in fact, goes to zero.

The last two parts of the remark follow from the first part and the proof of the first part immediately.
$\eop$

Of course, the same happens, if in Remark 3 we replace $X^{2d}_{k}$ by $X^{d}_{k}-X^{2d}_{k}$.

From now on, because of Remark 3, we will make the following assumption. It is not really needed in any of the proofs
but it guarantees that our approximations of Feynman path integrals (see below) indeed approximate them
in a reasonable sense without which
there is no point in what we do (keep in mind that Feynman path integral, especially the first version
that actually talks about paths, does not have mathematically rigorous definition). And in physics, everything tends to be
continuous (although sometimes we have poles or Dirac's delta kind of behaviours, which famously couse problems).

\medskip

{\bf Assumption 3}. We assume, that $f_{v}:\R\rightarrow\R$ is continuous.

\medskip

This assumption holds for the harmonic oscillator. 

And now to Feynman path integrals.
The idea here is a bit similar than the idea behind Bochner spaces.
Alternatively one can think that we work in an ultraproduct of spaces $H^{d}_{k}$ together with the field of
complex numbers as we did to some extend in [HH1] (one needs to be careful here, complex multiplication does not behave well in a metric ultraproduct).
We write $\C^{\o}$
for the ring of functions $f:\o\rightarrow\C$ with coordinatewise addition and multiplication.
The ring $\C^{\o}$ is not a field since $f\in\C^{\o}$ has the multiplicative inverse $f^{-1}$
iff for all $n<\o$, $f(n)\ne 0$. If $g,f\in\C^{\o}$ and $f^{-1}$ exists, we write as usually $g/f$
for $gf^{-1}$. Since ultimately we are only interested in what happens for almost all $n<\o$,
we could use this notation also when $f(n)\ne 0$
for almost all $n$. But there is no need for this.

Letting
$N_{\n}=\{ f\in\C^{\o}\vert\ \n\subseteq f\}$, $\n\in\C^{<\o}$, be the basic open sets, we get a metrizable topology
to $\C^{\o}$. We write $\R^{\o}$ for the subring of functions $f:\o\rightarrow\R$ and in $\R^{\o}$ the topology is that
induced from $\C^{\o}$. In $\C^{\o}$ we can also define other functions, e.g. $e^{x}$, coordinatewise
and the same for constants e.g. $i$ becomes the constant $i$ function etc. Some of our object are not
constant nor functions but we can still give them the meaning the obvious way, e.g.
$\Delta t$ becomes the function $k\mapsto t/n^{*}(k)$
and $N^{d}_{k}$ becomes the function $k\mapsto N^{d}_{k}$. Notice that we could have defined $e^{x}$
also by using the usual formula and got the same function. 

We say that $k<\o$ is $\d$-good if for all $p<k$, $\d (p)\le\d (k)$. We say that $\d$ is good if
the set of $\d$-good natural numbers is unbounded in $\o$. All functions $\d$ that we will look are good
and so for the rest of this section, we assume that $\d$ is good.

Let us now look what the set $P=P^{\d}$ of our paths is. We get this set as a union of two
sets $P^{2\d}$ and $P^{1\d}$
We let the set $P^{2\d}$ of paths to be the set of all approximations of continuously differentiable paths
$p:[0,t]\rightarrow\R$ from Remark 3.
If $p=((p^{k}_{j})_{j\le n^{*}(k)})_{k<\o}\in P$ is a path, we write $p^{k}=(p^{k}_{j})_{j\le n^{*}(k)}$.
Notice that if $k$ is $\d$-good,
$p^{k}$ determines $p^{r}$ for all $r<p$ ($p^{k}$ determines $p^{r}$ if  $\sqrt{N^{d}_{r}}$ divides $\sqrt{N^{d}_{k}}$).
Thus we have projections i.e.
I can write $Pr_{kr}(p^{k})$ for $p^{r}$ in the case $k$ is $\d$-good and $r<k$.
We add a topology to $P^{2\d}$ by letting $N_{x}=N_{x}^{k}$ be the basic open sets, where $k$ is $\d$-good and
$x\in X^{2d}_{k}\times (X^{2d}_{k})^{n^{*}(k)-1}\times X^{2d}_{k}$ and $p\in N_{x}$ if 
$p^{k}=(p^{k}_{j})_{j\le n^{*}(k)}=x$.
I will usually drop $k$ from this notation although one is not necessarily able to read it out from $x$
(often one can).
Notice that
basic open sets are clopen and that the topology is metrizable.
I write $m(N_{x}^{k})=1/k$.

The set $P^{1\d}$, its topology, projections and $m(N^{k}_{x})$ are defined exactly as in the case of $P^{2\d}$ except that
we replace $X^{2d}_{k}$ by $X^{d}_{k}-X^{2d}_{k}$ everywhere (so here the approximations $p=(p^{k})_{k<\o}$ of continuously
differentiable paths are such that $p^{k}\in (X^{d}_{k}-X^{2d}_{k})^{n^{*}+1}$).
Notice that from $x$ (and $k$) in $N^{k}_{x}$ one  can read out if
$N^{k}_{x}$ is a basic open set of $P^{2\d}$ or $P^{1\d}$ (from $x$ one can read out $\d (k)$ and then with $k$, also $N^{d}_{k}$).
The same is true when $x$ is an argument in the projections.

Finaly, we let $P^{\d}=P^{2\d}\cup P^{1\d}$ and we let the set of basic open sets of $P^{\d}$ be the union
of the sets of basic open sets of $P^{2\d}$ and $P^{1\d}$.
So now each continuously differentiable path $p:[0,t]\rightarrow\R$ has two approximations,
one in $P^{2\d}$ and one in $P^{1\d}$ and when we calculate the path integral both are included in the integration.
So in a sense, each path will be counted twice in our integration. Obviously this is not a problem even in a philosophical
sense, since in the definition of Faynman propagator a normalizing constant is used, which
can take care of things like this - if needed.

Let $P^{*}$ be the completion of $P$.
Then $P^{*}$ is compact. We will work in $P$ although our definitions would look more natural
in the compact $P^{*}$.
If $X$ is a finite set of basic open sets, then $m(X)=max\{ m(N)\vert\ N\in X\}$.
Following the terminology from Riemann integrals, I call this the mesh of $X$.

And now to the measure. Our measure is not really a measure, the values are not from $\R$ but keep also in mind
the non-existence of the real Lebesgue-style measure for differentiable paths that fits the bill, see [Pe]. Our measure is
more or less Lebesgue-style
(i.e. translations do not change the measure, it looks the
same everywhere, at least for almost all $k$ which is all we care, see below)
only it is not exactly a measure. But
it allows us to define integrals needed here. One may want to compare our measure to Wiener measure,
which is used by physicists, when their approach to Feynman path integral is that of imaginary time, see [Me].
I let $\mu_{k}$ be the counting measure on
$$X^{**d}_{k}=(X^{2d}_{k}\times (X^{2d}_{k})^{n^{*}(k)-1}\times X^{2d}_{k})\cup ((X^{d}_{k}-X^{2d}_{k})\times (X^{d}_{k}-X^{2d}_{k})^{n^{*}(k)-1}\times (X^{d}_{k}-X^{2d}_{k}))$$
normalized so
that
$$\mu_{k}(X^{**d}_{k})=(N^{d}_{k}/2)^{2}\mu^{k}(X^{*d}_{k}),$$
where $\mu^{k}$ and $X^{*d}_{k}$ are as defined few paragraphs before Remark 3. Notice that
$$\mu_{k}(\{ 0\}\times (X^{2d}_{k})^{n^{*}(k)-1}\times\{ 0\} )=\mu^{k}((X^{2d}_{k})^{n^{*}(k)-1}),$$
cf. the discussion just before Remark 3.
Next I define a measure-like function $\mu^{*}$ to $P$ (and $P^{*}$)
so that it gets its values from $\R^{\o}$. It will be enough to define $\mu^{*}$ for
basic open sets only, we will 'Riemann integrate' uniformly continuous functions.

So we want to define $\mu^{*}(N^{k}_{x})(r)$, $k$ $\d$-good and $r<\o$.
If  $r>k$, it does not matter how we define this because of the way we defined the
topology in $P$ and $\C^{\o}$. The natural choice is the following but see Remark 4 below:
We write
$$Y_{x}^{r}=Y^{kr}_{x}=\{ p^{r}\vert\ p\in N^{k}_{x}\} .$$
Then we let
$$\mu^{*}(N^{k}_{x})(r)=\mu_{r}(Y^{r}_{x} ).$$

Next suppose $r<k<\o$. List all $\d$-good $r<q\le k$ in an increasing order as $q_{i}$, $1\le i\le p$,
and let $q_{0}=r$. For all $i<p$, let
$$m_{x}^{i}=1/\vert Y^{q_{i+1}}_{Pr_{kq_{i}}(x)}\vert$$
and then
$$\mu^{*}(N_{x})(r)=(\Pi_{i=0}^{p-1}m_{x}^{i})\mu_{r}(\{ Pr_{kr}(x)\} ).$$
This is not the only way to define this, but we need something that guarantees that the paths are calculated
the right number of times i.e. with the right weight (cf. the last case below and see also Remark 4 below)
and this was the first method that came to my mind.
Finally, we are left with the case $r=k$:
$$(5)\ \ \ \ \ \mu^{*}(N_{x})(k)=\mu_{k}(\{ x\}).$$

{\bf Remark 4}. Since we did not care much how $\mu^{*}(N^{k}_{x})(r)$ was defined in the case $r>k$,
the following is possible. Let $X$ be a finite set of disjoint basic open sets  of $P$.
Then it is possible that
$$(6)\ \ \ \ \ \sum_{N\in X}\mu^{*}(N)(r)>\mu_{r}(X^{**d}_{r}).$$
This is because it is possible that there are
$q\in N^{k}_{x}\in X$ and $p\in N^{k}_{y}\in X$ such that $q_{r}=p_{r}$ but
$x\ne y$ and then  this $q_{r}=p_{r}$ is counted more than once in the sum.
However, for every $X$ this can happen only for finitely many $r$
(it does not happen when $N_{r}$ and thus $N^{d}_{r}$ is divisible enough)
and it can never happen for those $r$ such that
$1/r$ is greater or equal to the mesh of $X$. If fact for these $r$ more is true:

(7) If in addition we assume that $X$ is such that for all $p,q\in\cup X$,
$p^{r}=q^{r}=x$, then $\sum_{N\in X}\mu^{*}(N)(r)\le\mu_{r}(\{ x\} )$ always
and $\sum_{N\in X}\mu^{*}(N)(r)=\mu_{r}(\{ x\} )$ if and only if for all $p\in P$,
$p^{r}=x$ implies $p\in\cup X$.

\noindent
This makes our integrals work. In fact we can choose $\mu^{*}$ in any way as long as it satisfies
(5) and (7) and it will work in our calculations and give the same results.

Notice that from (6) it follows that it is possible that $X$ and $Y$
are finite sets of disjoint basic open sets, $\cup X=\cup Y$
but $\sum_{N\in X}\mu^{*}(N)\ne \sum_{N\in Y}\mu^{*}(N)$.
One can get rid of all cases of (6) if one wants: One can call $\mu^{*}$ a pre-measure like function and then
define a measure like function $\mu^{**}$ so that for all basic open sets $N$,
$\mu^{**}(N)=\int_{N}1dp$, where $\int_{N}$ is defined as below using $\mu^{*}$ and
$1$ is the function that for all $q\in P$, $1(q)(k)=1$ for all $k<\o$.
Then for unions of finite sets of disjoint basic open sets one can define the measure
as the sum of the measures of the basic open sets. And one can continue developing theory
for this $\mu^{**}$. But we will not need this, in fact,
in the definition of our integrals (in the cases in which we use it), it does not matter whether we use
$\mu^{*}$ or $\mu^{**}$ and thus for simplicity we use $\mu^{*}$.

Integration is defined as follows: Suppose $f:P\rightarrow\C^{\o}$ is uniformly continuous i.e.
for all $k$ there is $m$ such that if $N$ is a basic open set, $m(N)<1/m$ and $x,y\in N$,
then $f(x)(k)=f(y)(k)$. Now suppose that
$P_{0}\subseteq P$ is a finite union of basic open sets (in compact $P^{*}$ this makes more sense) and
$X$ is a finite partition of $P_{0}$ consisting of basic open sets
(i.e. $\cup X=P_{0}$, every element of $X$ is a basic open set and if $N,N'\in X$ and $N\ne N'$, then $N\cap N'=\empty$).
For all $N\in X$, pick $p_{N}\in N$.
In fact, we can fix these choices $p_{N}\in N$ for all basec open sets $N$.
Then we write $I_{f}(X)=\sum_{N\in X}f(p_{N})\mu^{*}(N)$ (notice that here we calculate in $\C^{\o}$).
By $\int_{P_{0}}fdp$ we mean the limit of $I_{f}(X)$, when the mesh of $X$ goes to zero (limit in the topology
we defined for $\C^{\o}$ above).
The assumption that $P_{0}$ is a finite union of basic open sets guarantees the existence
of these partitions with arbitrarily small mesh. Thus
it is easy to see that
this limit exists and it does not depend on the choice of the partitions $X$ (as long as their mesh goes to zero)
nor on the choices of the paths $p_{N}$ (since $f$ was assumed to be uniformly continuous).
For sets that are not finite unions of basic open sets,
things get trickier. We will not need these cases in this paper,
we will integrate only over sets that are finite unions of basic open sets.

Notice that for all $c\in\C^{\o}$, uniformly continuous $f$ and $g$ and finite unions $P_{0}$ and $P_{1}$ of basic open sets
such that $P_{0}\cap P_{1}=\empty$,
$c\int_{P_{0}}fdp=\int_{P_{0}}cfdp$, $\int_{P_{0}}(f+g)dp=(\int_{P_{0}}fdp)+(\int_{P_{0}}gdp)$ and
$\int_{P_{0}\cup P_{1}}fdp=(\int_{P_{0}}fdp)+(\int_{P_{1}}fdp)$.

Next we define action $S(p)$ of a path $p\in P$. Obviously
$$S(p)(k)=\Delta t\sum_{j=1}^{n^{*}(k)}S^{*}(p^{k}_{j-1},p^{k}_{j},\Delta t),$$
where $\Delta t=t/n^{*}(k)$.
This depends on $\d$ and we will use this definition for various $\d$'s. Our final $\d =\d^{*}$
will be defined in Subsection 3.2.
Notice that $S$ is uniformly continuous and thus also the function
$(e^{iS(p)/\hbar})(k)=e^{i((S(p))(k))/\hbar}$ is uniformly continuous.

Finally, in Theorem 1, also an ultrafilter $D^{*}$ is needed.
In Subsection 3.2 we will define a function $\tau :\o\rightarrow\o$ (which is closely connected to
$\delta^{*}$) and then
we let $D^{*}$ be any ultrafilter on $\o$ such that for all $n$,
$$\{ k<\o\vert\ \tau(k)\ge n\}\in D^{*}.$$

Now we are ready to state the main theorem except that we need one more assumption.
We will be able to calculate the propagator much as in [HH1], except that we can not do the cleaning of the formula
(and thus justify the final limit process) that we were able to do
in the end in [HH1]. This is because it is very hard to argue that the big sum of the formulas
$$e^{i\Delta t(\sum_{j=1}^{n^{*}}S^{*}(x_{j-1},x_{j},\Delta t))/\hbar}$$
that we will end up studying in spaces $H^{d}_{k}$ behaves well when $n$ goes to infinity,
When $n$ increases,
the sum may start oscillating more and more wildly, so fast that I cannot handle it and I cannot compensate this e.g. by looking
it in smaller and smaller intervals (because Trotter-Kato does not allow this, it gives only a pointwise convergence).
In [HH1] we did not have this problem, everything was nicely continuous (modular the spiking effect that was
easy to handle). This oscillation problem is also the reason why it is very hard to get rid of
the $lim_{r\rightarrow\infty}$ part from the formula  in Theorem 1
and $lim_{p\rightarrow\infty}$ part from the formula in the end of Section 1.1. We will return to this in Section 3.1.
Thus we will need something that replaces
the (nice) continuity of the formula from [HH1].
Our assumptions appears to be natural in our context that arises from physics where functions tend to be continuous
(well, as already mentioned, they may have poles or Dirac's delta kind of behaviour). The assumption is true e.g.
for the harmonic oscillator in the non-trivial cases
i.e. when the time evolution operator is not
the identity operator (if it is the identity operator, it does not have a kernel).

\medskip

{\bf Assumption 4}: We assume that (the Hamiltonian is such that) at our fixed time $t>0$,
the kernel $K(x_{0},x_{1},t)$ of the time evolution operator $K^{t}$ in $L_{2}(\R )$ is continuous
as a function $\R^{2}\rightarrow\C$ (in particular we assume that the
time evolution operator has a kernel).

\medskip

It is easy to see that our calculations work also in many cases in which this Assumption 4 fails but
without any additional assumptions, I do not see how to push the calculation through.
However, when we calculate $K(x_{0},x_{1},t)$ for specific $x_{0}$ and $x_{1}$,
to make our technique work, it is enough to assume that
the operator has this kernel and that it is continuous at $(x_{0},x_{1})$.
Also in Theorem 1 below, keep in mind our requirements for the choice of units
i.e. units are chosen so that the requirements hold for time $t$ and the mass $m$ of the particle
including the requirement that $th/2m =1$. For other integer values for $th/2m$, one needs to make changes to
many places. E.g. we need to divide $X^{d}_{k}$ into $th/m$ many pieces and not just to $2$ (kind of modulo $th/m$)
and thus also $P^{\d}$ becomes a union of $th/m$ many subspaces. Also some changes are need in the
definitions of our measures and the measure like function etc. I have not checked all the details
but the cases $th/2m >1$ are more complicated than our $th/2m=1$ case, which is the reason why
we chose the scaling the way we did.

So in Theorem 1 below, we use our original units i.e. units from the spaces $H_{k}$ as we did in Remark 3 above. In particular,
$(m/2\pi\hbar t)^{1/2}$ is calculated in the original units.
Also in the definition of $P_{y_{0},y_{1},r}$ these units are used (we will use them also in $X^{d}_{k}$).
In fact, in the proof of the theorem, we will switch to use entirely these units
in the formulas (but we will not make changes to the spaces $H^{d}_{k}$).
This does not affect to the numerical value of
$S(p)/\hbar$ since this is dimensionless. So if one wants, one can calculate $\int_{P_{y_{0},y_{1},r}}e^{iS(p)/\hbar}dp$ in the spaces
$H^{d}_{k}$ using units from there.

\th 1 Theorem. Under the assumptions 1-4, there are (good) $\d^{*}:\o\rightarrow\o$ and
an ultrafilter $D^{*}$ on $\o$ such that
for all rational numbers $y_{0}$ and $y_{1}$, the kernel $K(y_{0},y_{1},t)$ is
$$(m/2\pi\hbar t)^{1/2}lim_{r\rightarrow\infty}lim_{D^{*}}(2^{-2\d^{*}(k)}r^{2}/4\sqrt{N_{k}})\int_{P_{y_{0},y_{1},r}}e^{iS(p)/\hbar}dp,$$
where $r$ is a non-zero natural number and $P_{y_{0},y_{1},r}$ is the set of all paths $p\in P$
such that $y_{0}-(1/r)\le p^{k}_{0}<y_{0}+(1/r)$ and $y_{1}-(1/r)\le p^{k}_{n^{*}(k)}<y_{1}+(1/r)$
for all large enough $k$.

Notice that in the definition of $P_{y_{0},y_{1},r}$ above, what is large enough $k$ is
determined by $y_{0}$, $y_{1}$ and $r$ 
i.e. it does not depend on $p$ and thus $P_{y_{0},y_{1},r}$ is a finite union of basic open sets.
The index $k$ is (certainly) large enough if $k>1$ and $\vert \sqrt{N_{k}}y_{0}\vert$, $\vert \sqrt{N_{k}}y_{1}\vert$
and $\sqrt{N_{k}}/r$ are natural numbers and $<(N_{k}/4)-1$.
Also keep in mind that not only $\int_{P_{y_{0},y_{1},r}}e^{iS(p)/\hbar}dp$ but also $2^{-2\d^{*}(k)}r^{2}/4\sqrt{N_{k}}$
is a  function $\o\rightarrow\C$.

\chapter{3 The proof of the main theorem}

This section is the proof of the main theorem and thus we assume that Assumptions 1-4 hold
(or at least Assumptions 1, 2 and 4).
The proof is a straight forward calculation.
In the end of Section 3.1, I will return to the $lim_{p\rightarrow\infty}$
problem that we discussed in the end of Section 1.1.
We start with the case when $\d$ is a constant function.

\chapter{3.1 The constant $\d$ case}

We pick (large) $s<\o$ and excluding the very end of this subsection,
we let $\d=\d_{s}$ be the constant function $\d_{s}(n)=s$ for all $n<\o$.
In the very end of this subsection we let $s$ run to infinity.

I start by calculating the propagators $<x_{0}\vert L^{t}_{k}\vert x_{1}>$ for $Q$-eigenvectors $x_{0}>$ and $x_{1}>$ in $H^{d}_{k}$.
We look first the case in which their eigenvalues are in $X^{2d}_{k}$. Recall that $X^{2d}_{k}$ is the set of all $x\in X^{d}_{k}$ such that
$x/2\in X^{d}_{k}$.
This is easy, using Lemma 1, the definition of $e^{-itV_{N^{d}_{k}}/\hbar 2^{2\d (k)}}$ i.e. that
$$e^{-itV_{N^{d}_{k}}/\hbar 2^{2\d (k)}}(y>)=e^{-i\Delta tf^{d}_{v}(y)/\hbar}\vert y>$$ 
and recalling
$S^{*}(x,y,t)=(m(y-x)^{2}/2t^{2})-f_{v}^{d}(x)$, $n^{*}=n^{*}(k)=2^{2\d (k)}=2^{2s}$ and $\Delta t=t/2^{2\d (k)}=t/2^{2s}$, we get the following:

\th 3 Lemma.
The propagator $<x_{n^{*}}\vert L^{t}_{k}\vert x_{0}>$, $x_{0},x_{n^{*}}\in X^{2d}_{k}$, is
$$((N_{k}^{d})^{-1/2}(\Delta th/m)(m/2\pi i\hbar\Delta t)^{1/2})^{n^{*}}(\sum_{x_{1}\in X^{2d}_{k}}...$$
$$...\sum_{x_{n^{*}-1}\in X^{2d}_{k}}e^{i\Delta t(\sum_{j=1}^{n^{*}}S^{*}(x_{j-1},x_{j},\Delta t))/\hbar}).$$

\proof
We can see this as follows: Clearly, by Lemma 1 and the definition of $e^{-itV_{N^{d}_{k}}/\hbar 2^{2\d (k)}}$,
$$L^{\Delta t}_{k}\vert x_{0}>=(N_{k}^{d})^{-1/2}(\Delta th/m)(m/2\pi i\hbar\Delta t)^{1/2}(\sum_{x_{1}\in X^{2d}_{k}}e^{i\Delta t S^{*}(x_{0},x_{1},\Delta t))/\hbar}\vert x_{1}>).$$
Here we can restrict the summation to the elements $x_{1}\in X^{2d}_{k}$
because if $x_{1}\in X^{d}_{k}-X^{2d}_{k}$, then $<x_{1}\vert L^{\Delta t}_{k}\vert x_{0}>=0$ by Lemma 1.
By repeating,
$$L^{2\Delta t}_{k}\vert x_{0}>=((N_{k}^{d})^{-1/2}(\Delta th/m)(m/2\pi i\hbar\Delta t)^{1/2})^{2}(\sum_{x_{1}\in X^{2d}_{k}}$$
$$\sum_{x_{2}\in X^{2d}_{k}}e^{i\Delta t (S^{*}(x_{0},x_{1},\Delta t)+S^{*}(x_{1},x_{2},\Delta t))/\hbar}\vert x_{2}>).$$
Here keep in mind 'Fubini' i.e. that
$$\sum_{x_{1}\in X^{2d}_{k}}\sum_{x_{2}\in X^{2d}_{k}}=\sum_{(x_{1},x_{2})\in (X^{2d}_{k})^{2}}=\sum_{x_{2}\in X^{2d}_{k}}\sum_{x_{1}\in X^{2d}_{k}}.$$
We will use this also later.
And so eventually we get the formula in the lemma for $<x_{n^{*}}\vert L^{t}_{k}\vert x_{0}>$ by noticing that
if $z_{i}$, $i<N^{d}_{k}$ enumerate $X^{d}_{k}$, then
for $i<N^{d}_{k}$,
$<z_{i}\vert (\sum_{j <N^{d}_{k}}a_{j}\vert z_{j}>)>=a_{i}$. $\eop$

If $x_{0},x_{n^{*}}\in X^{d}_{k}-X^{2d}_{k}$, then we get the propagator from the formula above by replacing $X^{2d}_{k}$ by $X^{d}_{k}-X^{2d}_{k}$ everywhere
and we also notice the following:

\th 4 Lemma.
If one of $x_{0}$ and $x_{n^{*}}$ is in $X^{2d}_{k}$ and the other one is not, then the propagator $<x_{n^{*}}\vert L^{t}_{k}\vert x_{0}>$ is $0$.

\proof Immediate by Lemma 1 and the proof of Lemma 3 (in the products
$$\Pi_{i=1}^{n^{*}}<x_{i-1}\vert L^{\Delta t}_{k}\vert x_{i}>,$$
by Lemma 1, at least one of the values is $0$ independent of the choice of the path $(x_{i})_{0<i<n^{*}}$). $\eop$

Let us first take a closer look at the formula from Lemma 3.
First let us write the formula in the 'original units' of $H_{k}$
i.e. we write $h$ in the units of $H_{k}$ and thus replace $h$ by $2^{2\d (k)}h$ in the formula in the part that is before the sums,
we do not need to touch to the formula after the sums since that part i.e. the formula in the exponent is dimensionless and thus nothing happens.
We also use the fact that $\Delta t2^{2\d (k)}h=th$.
Notice that
this does not change the numerical value
of the propagator. From now on we will use the original units of $H_{k}$ so e.g. the new $X^{d}_{k}$ is
$\{ 2^{-\d (k)}x\vert\ x\in X^{d}_{k}\}$.
With all this the formula from Lemma 3 gets the form
$$((N_{k}^{d})^{-1/2}(th/m)(m/2\pi i\hbar t)^{1/2})^{n^{*}}(\sum_{x_{1}\in X^{2d}_{k}}...$$
$$...\sum_{x_{n^{*}-1}\in X^{2d}_{k}}e^{i\Delta t(\sum_{j=1}^{n^{*}}S^{*}(x_{j-1},x_{j},\Delta t))/\hbar}).$$

Then we will
use integral notation for the sums over $X^{2d}_{k}$.
Recall that we defined the measure so that
$dx_{1}...x_{n^{*}-1}=\sqrt{2}(N^{d}_{k}/2)^{-n^{*}/2}$.
Then the propagator $<x_{n^{*}}\vert L^{t}_{k}\vert x_{0}>$ becomes, assuming $x_{0},x_{n^{*}}\in X^{2d}_{k}$
(keep in mind that $th/m=2$ and $2\pi\hbar =h$),
$$(8)\ \ \ \ \ (m/2\pi \hbar t)^{1/2}\int_{(X^{2d}_{k})^{n^{*}-1}}e^{i\Delta t(\sum_{j=1}^{n^{*}}S^{*}(x_{j-1},x_{j},\Delta t))/\hbar}dx_{1}...x_{n^{*}-1},$$
where we have assumed that $n^{*}$ is divisible by $8$ which is almost always true
also in the  case when $\d$ is got  by diagonalizing (where in the formula $i$ ends up, depends on how we
split the time and our choice makes it disappear altogether from this formula - but it still lurks somewhere).
We have not pushed the term $(m/2\pi \hbar t)^{1/2}$ into the measure for two reasons: First of all,
physicists seem to use often this very term as the normalizing constant $A(t)$
in their calculations (sometimes multiplied with an infinite constant) and secondly,
although it is very difficult here to track the dimensions, one expects that
$$\int_{(X^{2d}_{k})^{n^{*}-1}}e^{i\Delta t(\sum_{j=1}^{n^{*}}S^{*}(x_{j-1},x_{j},\Delta t))/\hbar}dx_{1}...x_{n^{*}-1}$$
is dimensionless (since the propagator should have dimension $1/length$) and so this term $(m/2\pi \hbar t)^{1/2}$ carries the unit
of the propagator naturally.
I find it an interesting question to figure out
the role and behaviour of the units in these calculations (including the calculations physicists do).
The case when $x_{0},x_{n^{*}}\in X^{d}_{k}-X^{2d}_{k}$ is similar.

Let us go back to units. Above we started using in formulas units from $H_{k}$. However,
still because we do the calculations in $H^{d}_{k}$ instead of $H_{k}$ the value may be wrong.
Indeed, since the unit of the propagator is $1/length$, one expects that the numerical value
one gets in $H^{d}_{k}$ is $2^{d(k)}$ times too small. There are also other problems which we will handle as we
did in [HH1].
Since all this is rather delicate,
let us look at the details. Keep in mind that our $\d$ is the constant  function $\d (k)=s$ for all $k$.
As before we let $L_{2}(\R )$ be the standard model that uses units from $H_{k}$.
We let $L^{s}_{2}(\R )$ be the standard model we get by using units we originally used in $H^{d}_{k}$,
in particular, $h$ in $L^{s}_{2}(\R )$ is  $2^{2s}$ times what it is in $L_{2}(\R )$
and if we write $f^{s}_{v}$ for $f_{v}$ in $L^{s}_{2}(\R )$, then
$f^{s}_{v}(x)=2^{2s}f_{v}(2^{-s}x)$. We define an isometric isomorphism $R=R^{s}$ from $L_{2}(\R )$
to $L^{s}_{2}(\R )$ so that $R(\phi )(x)=2^{-s/2}\phi (2^{-s}x)$. Now if we write
$V^{s}$ for the operator in $L^{s}_{2}(\R )$ that we get from $f^{s}_{v}$,
$$e^{-itV^{s}/\hbar}=R\circ e^{-itV/\hbar}\circ R^{-1}$$
(for all $t$ and keep in mind that the numerical value of  $\hbar$ in $L^{s}_{2}(\R )$ is $2^{2s}$ times the numerical
value of $\hbar$ in  $L_{2}(\R )$). 
Also if $W$ is an operator in $L_{2}(\R )$ with kernel $K_{W}(x,y)$,
then $R\circ W\circ R^{-1}$ is an operator in $L^{s}_{2}(\R )$ with kernel
$K'(x,y)=2^{-s}K_{W}(2^{-s}x,2^{-s}y)$. Now if one looks what the kernel for the time evolution operator
for the free particle is (see Lemma 1 and the paragraph just before it),
one notices that $G$ preserves also the time evolution operator for the free particle
(for all $t$).
Now let us
write  $L^{*t}_{s}$ for the operator
$$(e^{-itP^{2}/2^{2s}2m\hbar}e^{-itV/2^{2s}\hbar})^{2^{2s}}$$
in $L_{2}(\R )$ and $L^{s*t}_{s}$ for the same operator in $L^{s}_{2}(\R )$.
Then $L^{s*t}_{s}=R\circ L^{*t}_{s}\circ R^{-1}$. We recall that if we calculate
the kernel of $L^{s*t}_{s}$ using approximations $H^{d}_{k}$, then to get the kernel
of $L^{*t}_{s}$, we need to multiply the former by $2^{s}$. And the same is true for other kernels.
However we will need this only in one remark below
This is because of the order in
which we will take our limits.
Below, the spaces $L^{s}_{2}(\R )$ and the isometries $R^{s}$, $s<\o$, will be needed frequently.

Now we pick an ultrafilter $D^{s}$ such that for all $n<\o$,
$\{ N^{d}_{k}\vert\ k>n\}\in D^{s}$ and notice that
the spaces $H^{d}_{k}$ approximate, when in them we use the original units of $H_{k}^{d}$, the space $L^{s}_{2}(\R )$ as
the spaces $H_{N}$ approximated the space $L_{2}(\R )$ (they have the same units and $D^{s}$ satisfies the
requirements for $D$), see Sections 1.2 and 2.
We write $F^{s}$ for the embedding of $L^{s}_{2}(\R )$
to the metric ultraproduct of the spaces $H_{k}^{d}$ (taken using $D^{s}$) and we write
$H^{s}_{u}$ for the image of $F^{s}$. By $F^{s}_{k}$ we mean the functions from which the function $F^{s}$
is build in Section 1.2.

We write  $L^{s*t}_{s}$ also for the
image of $L^{s*t}_{s}$ under the embedding $F^{s}$ in $H^{s}_{u}$.
Then by Assumption 2, we can use these $L^{s*t}_{s}$ to calculate the propagator $K(y_{0},y_{1},t)$
in $L_{2}(\R )$. And by Remark 2, our approximations tell what happens with $L^{s*t}_{s}$.
The calculation goes much as in [HH1]: 
We write $\phi_{y,p}\in L_{2}(\R )$, $p>0$ a natural number,
for the function such that $\phi_{y,p}(x)=1$ if $x\in [y-(1/p),y+(1/p)[$ and otherwise $\phi_{y,p}(x)=0$ (or a test function that approximates
this function better and better, see [HH1] or Subsection 1.2). We write $\phi^{s}_{y,p}$ for
$R^{s}(\phi^{s}_{y,p} )$.

Then by Assumption 4, since
$$<\phi_{y_{1},p}\vert K^{t}\vert\phi_{y_{0},p}>=\int_{[y_{1}-(1/p),y_{1}+(1/p)[}\int_{[y_{0}-(1/p),y_{0}+(1/p)[}K(x,y,t)dxdy,$$
we get that
$$(9)\ \ \ \ \ K(y_{0},y_{1},t)=lim_{p\rightarrow\infty}(p/2)^{2}<\phi_{y_{1},p}\vert K^{t}\vert\phi_{y_{0},p}>,$$
where $(p/2)^{2}$ comes from the fact that when we calculated $<\phi_{y_{1},p}\vert K^{t}\vert\phi_{y_{0},p}>$
above, we integrated twice over an interval of length $2/p$.
But now it follows from Assumption 2 and the fact that $R$ is an isometric isomorphism, that
$$(10)\ \ \ \ \ K(y_{0},y_{1},t)=lim_{p\rightarrow\infty}(p/2)^{2}lim_{s\rightarrow\infty}<\phi_{y_{1},p}\vert L^{*t}_{s}\vert\phi_{y_{0},p}>=$$
$$lim_{p\rightarrow\infty}(p/2)^{2}lim_{s\rightarrow\infty}<\phi^{s}_{y_{1},p}\vert L^{s*t}_{s}\vert\phi^{s}_{y_{0},p}>,$$
since $K^{t}(\phi_{y_{0},p})=lim_{s\rightarrow\infty}L^{*t}_{s}(\phi_{y_{0},p})$ and the isometry $R$ maps $L_{s}^{*t}$ to $L_{s}^{s*t}$.
We are left to calculate $<\phi^{s}_{y_{1},p}\vert L^{s*t}_{s}\vert\phi^{s}_{y_{0},p}>$ using our approximations.

Now
$$<\phi^{s}_{y_{1},p}\vert L^{s*t}_{s}\vert\phi^{s}_{y_{0},p}>=lim_{D^{s}}<F^{s}_{k}(\phi^{s}_{y_{1},p})\vert L^{t}_{k}\vert F^{s}_{k}(\phi^{s}_{y_{0},p})>,$$
where 
we think the expression $<F^{s}_{k}(\phi^{s}_{y_{1},p})\vert L^{t}_{k}\vert F^{s}_{k}(\phi^{s}_{y_{0},p})>$ as a function
$$k\mapsto <F^{s}_{k}(\phi^{s}_{y_{1},p})\vert L^{t}_{k}\vert F^{s}_{k}(\phi^{s}_{y_{0},p})>.$$
Below, similar expressions are also considered as functions.
We are left to put together the results from this subsection:

For rational numbers $x$ and natural numbers $p>0$, we write $X^{2d}_{k}(x,p)=\{ y\in X^{2d}_{k}\vert\ x-(1/p)\le y<x+(1/p)\}$,
keep in mind that in $H^{d}_{k}$ we are now using units from $H_{k}$.
Also the formula from Lemma 3 in the integral form is written using units from $H_{k}$, see the discussion after Lemma 4.
Since this did not change the numerical value of the formula, the changed version can be used here.
All this, including the formulas, is true also for $X^{d}_{k}-X^{2d}_{k}$ in place of $X^{2d}_{k}$
and we write $\overline{X}^{2d}_{k}(x,p)=\{ y\in X^{d}_{k}-X^{2d}_{k}\vert\ x-(1/p)\le y<x+(1/p)\}$
for rational numbers $x$ and natural numbers $p>0$. Also keep in mind that $x\in X^{2d}_{k}$ for $k$ large enough.

We have proved:
For rational numbers $y_{0}$ and $y_{1}$ and $k$ large enough
the inner product
$$<F^{s}_{k}(\phi^{s}_{y_{1},p})\vert L^{t}_{k}\vert F^{s}_{k}(\phi^{s}_{y_{0},p})>$$
is
$$(m/2\pi \hbar t)^{1/2}(2^{-s}/\sqrt{N^{d}_{k}})$$
$$((\sum_{x_{0}\in X^{2d}_{k}(y_{0},p)}\sum_{x_{n^{*}}\in X^{2d}_{k}(y_{1},p)}\int_{(X_{k}^{2d})^{n^{*}-1}}e^{i\Delta t(\sum_{j=1}^{n^{*}}S^{*}(x_{j-1},x_{j},\Delta t)/\hbar)}dx_{1}...x_{n^{*}-1})+$$
$$(\sum_{x_{0}\in \overline{X}^{2d}_{k}(y_{0},p)}\sum_{x_{n^{*}}\in \overline{X}^{2d}_{k}(y_{1},p)}\int_{(X^{d}_{k}-X_{k}^{2d})^{n^{*}-1}}e^{i\Delta t(\sum_{j=1}^{n^{*}}S^{*}(x_{j-1},x_{j},\Delta t)/\hbar)}dx_{1}...x_{n^{*}-1}))$$
In the formula above, $1/\sqrt{N^{d}_{k}}$ comes from (the definition of) the function $F^{s}_{k}$
and $2^{-s}$ comes from the use of $R^{s}$ in the definition of $\phi^{s}_{y,p}$ and the rest comes from (8) above and
the last sentence of the paragraph from which (8) can be found.

We can rewrite this formula by letting $X^{**d}_{k}(x,y,p)$ be
$$\{ (x_{0},...,x_{n^{*}})\in X^{**d}_{k}\vert\ x-(1/p)\le x_{0}<x+(1/p),\ y-(1/p)\le x_{n^{*}}<y+(1/p)\}$$
for rational numbers $x$ and $y$ and natural numbers $p>0$.

\th 5 Lemma. For rational numbers $y_{0}$ and $y_{1}$ and $k$ large enough
the inner product
$$<F^{s}_{k}(\phi^{s}_{y_{1},p})\vert L^{t}_{k}\vert F^{s}_{k}(\phi^{s}_{y_{0},p})>$$
is
$$(m/2\pi \hbar t)^{1/2}(2^{-s}/\sqrt{N^{d}_{k}})$$
$$\int_{X_{k}^{**d}(y_{0},y_{1},p)}e^{i\Delta t(\sum_{j=1}^{n^{*}}S^{*}(x_{j-1},x_{j},\Delta t)/\hbar)}dx_{0}...x_{n^{*}}$$
and thus
the propagator
$K(y_{0},y_{1},t)$ in $L_{2}(\R )$ is
$$(11)\ \ \ \ \ (m/2\pi \hbar t)^{1/2}lim_{p\rightarrow\infty}lim_{s\rightarrow\infty}lim_{D^{*}}(2^{-2s}p^{2}/4\sqrt{N_{k}})$$
$$\int_{X_{k}^{**s}(y_{0},y_{1},p)}e^{i\Delta t(\sum_{j=1}^{2^{2s}}S^{*}(x_{j-1},x_{j},\Delta t)/\hbar)}dx_{0}...x_{2^{2s}},$$
where, for clarification, we have written $X^{**s}_{k}(y_{0},y_{1},p)$ for
the set $X^{**d}_{k}(y_{0},y_{1},p)$ in the case when $\d$ is such that $\d (k)=s$ (this makes the role of $s$ explicit).
$\eop$

In the latter formula, $p^{2}/4$ comes from $(9)$ above and notice that $2^{-s}/\sqrt{N^{d}_{k}}=2^{-2s}/\sqrt{N_{k}}$.
Notice also that for $k$ large enough, $y_{0},y_{1}\in X^{2d}_{k}$ (since they are rational numbers).

Formula (11) above looks a bit like the  limit of (1) when $n$ goes to infinity (by Fubini, which, as pointed out above, holds for our
discrete integrals). But we would like to massage it to a form that looks even more like the limit of (1): 
If we can  approximate $e^{itH/\hbar}$ very well, $lim_{D^{s}}$ in the formula above can be replaced by $lim_{k\rightarrow\infty}$ but
I do not see how to do this in general. Also in some cases one can get rid of $lim_{p\rightarrow\infty}$ but again I do not see how to do this
in general, see the discussion just before Assumption 4.
But in the case of the free particle (i.e. $f_{v}\equiv 0$)
for rational numbers $y_{0}$ and $y_{1}$,
the formula for $K(y_{0},y_{1},t)$ reduces to the following form
$$(12)\ \ \ \ \ \ \ \ \ \ \ \ \ (m/2\pi \hbar t)^{1/2}lim_{s\rightarrow\infty}\ \ \ \ \ \ \ \ \ \ \ \ \ \ \ $$
$$lim_{k\rightarrow\infty}(2^{2s}\sqrt{N_{k}}/2)\int_{(X^{2s}_{k})^{n^{*}(k)-1}}e^{i\Delta t(\sum_{j=1}^{2^{2s}}S^{*}(x_{j-1},x_{j},\Delta t)/\hbar)}dx_{1}...x_{2^{2s}-1},$$
where $x_{0}=y_{0}$ and $x_{n^{*}(k)}=y_{1}$.
However, the reason why this holds (see below) is the very simple nature of this case i.e. the approximations are 'too good'.
However, I do not see any reason why this could not happen in other cases as well
(but the reason must be more complicated), thus I ask:

{\bf Question 2}. Does something like (12) give the propagator in some cases other than the free  particle?

One way of seeing (12) is the following: In [HH1] we showed that the propagator $K(y_{0},y_{1},t)$ is
$$lim_{k\rightarrow\infty}(\sqrt{N_{k}}/2)<y_{1}\vert K^{t}_{k}\vert y_{0}>$$
(again here $2$ is $th/m$ and one can see this also simply from Lemma 1 and the paragraph just before it,
in fact, one could replace this limit by 'large enough $k$', since
$y_{0}$ and $y_{1}$ are assumed to be rational numbers, but this would be unreasonable in cases
other than this free particle case). But
the same is true for $L^{s}_{2}(\R )$ and thus
one can replace $N_{k}$ by $N^{d}_{k}$
and $K^{t}_{k}$ by $K^{t}_{N^{d}_{k}}$ (assuming $\d (k)=s$)
in this free particle case and get the propagator $K^{s}(y_{0},y_{1},t)$ in $L^{s}_{2}(\R )$.
Above we sow that
$$K(y_{0},y_{1},t)=2^{s}K^{s}(2^{s}y_{0},2^{s}y_{1},t).$$
Now in $H^{d}_{k}$ we have switch to use our original units and
e.g. $2^{s}y_{0}$ is $y_{0}$ in $X^{d}_{k}$.
Also $2^{s}\sqrt{N^{s}_{k}}=2^{2s}\sqrt{N_{k}}$
and clearly,
$K^{t}_{N^{d}_{k}}=L^{t}_{k}$ and thus we get
$$K(y_{0},y_{1},t)=lim_{k\rightarrow\infty}(2^{2s}\sqrt{N_{k}}/2)<y_{1}\vert L^{t}_{ks}\vert y_{0}>.$$
Finally, to get (12),
first we replace $<y_{1}\vert L^{t}_{ks}\vert y_{0}>$ by the formula we calculated for it above,
see (8) and again keep in mind that $y_{0}$ and $y_{1}$ belong to
$X^{2d}_{k}$ for all large enough $k$, and then we add the void  $lim_{s\rightarrow\infty}$. It is void because the sequence is  
a constant sequence in this free particle case,
in other cases it would not be a constant sequence.

\chapter{3.2 Diagonalization}

We can write the formula (11) above also in a form that looks a bit like
the first formulation (0) of Feynman path integral in the introduction above and this will be our last
observation i.e. 
we are ready to finish the proof of Theorem 1:
For this we redefine $\d =\d^{*}:\o\rightarrow\o$ by diagonalizing as follows:
In addition to $\d^{*}$,
we define also a function $\tau :\o\rightarrow\o$ and
we enumerate all functions $\psi_{y,r}$, $y$ a rational number and $r$ a non-zero natural number, from the previous subsection as
$\psi_{i}$, $i<\o$.
We write $\psi^{s}_{i}$ for $R^{s}(\psi_{i})$
Then
for all $k<\o$, we let $\tau (k)$ be the largest $0<s\le k$ for which there is $s\le s'\le k$ such that

(i) for all $i\le s$, in $L_{2}(\R )$, $\vert K^{t}(\psi_{i})-L^{*t}_{s'}(\psi_{i})\vert <1/s$,

(ii) for all $i,j\le s$,
$$\vert <\psi^{s'}_{i}\vert L^{s'*t}_{s'}\vert\psi^{s'}_{j}>-<F^{s'}_{k}(\psi^{s'}_{i})\vert L^{t}_{ks'}\vert F^{s'}_{k}(\psi^{s'}_{j})>\vert <1/s,$$

\noindent
(for $L^{t}_{ks'}$, see the paragraph just before Remark 2 and for
$L^{*t}_{s'}$ and $L^{s'*t}_{s'}$ see the previous subsection)
if there is such $s$
and then we let $\d^{*}(k)$ be the least $s'\ge s$ that satisfies (i) and (ii) for this $s$
and if there is no such $s$, we let
$\tau (k)=\d^{*}(k)=0$.
The key observation here is that for all
$n<\o$, $\{ k<\o\vert\ \tau (k)\ge n\}$ is unbounded
(and then also $\{ k<\o\vert\ \d^{*}(k)\ge n\}$ is unbounded and thus $\d^{*}$ is good):
Let $n<\o$ be given. Then by Assumption 2 there is $s'\ge n$
such that for all $i\le n$, in $L_{2}(\R )$, $\vert K^{t}(\psi_{i})-L^{*t}_{s'}(\psi_{i})\vert <1/n$.
But
then by Remark 2 and \L os's theorem,
the set of all $k$ such that for all $i,j\le n$,
$$\vert <\psi^{s'}_{i}\vert L^{s'*t}_{s'}\vert\psi^{s'}_{j}>-<F^{s'}_{k}(\psi^{s'}_{i})\vert L^{t}_{ks'}\vert F^{s'}_{k}(\psi^{s'}_{j})>\vert <1/n$$
belongs to $D^{s'}$.
Unboundedly many of these are $\ge s'$. Thus there are unboundedly many $k$ such that $\tau (k)\ge n$.
It follows that we can find an ultrafilter $D^{*}$ such that for all $n<\o$,
$\{ k<\o\vert\ \tau (k)>n\}\in D^{*}$.
Notice that since $R$ is an isometric isomorphism and $\vert\psi_{i}\vert_{2}\le 2$,
if $\tau (k)\ge n$, then for all $i,j\le n$,
$$\vert <\psi_{i}\vert K^{t}\vert\psi_{j}>-<F^{\d^{*}(k)}_{k}(\psi^{\d^{*}(k)}_{i})\vert L^{t}_{k\d^{*}(k)}\vert F^{\d^{*}(k)}_{k}(\psi^{\d^{*}(k)}_{j})>\vert <3/n.$$

Thus for all $\psi_{y_{0},r}$ and $\psi_{y_{1},r}$
$$(13)\ \ \ \ \ <\psi_{y_{1},r}\vert K^{t}\vert\psi_{y_{0},r}>=lim_{D^{*}}<F^{\d^{*}(k)}_{k}(\psi^{\d^{*}(k)}_{q_{1},r})\vert L^{t}_{k\d^{*}(k)}\vert F^{\d^{*}(k)}_{k}(\psi^{\d^{*}(k)}_{q_{0},r})>,$$
keep in mind that we think $<F^{\d^{*}(k)}_{k}(\psi^{\d^{*}(k)}_{q_{1},r})\vert L^{t}_{k\d^{*}(k)}\vert F^{\d^{*}(k)}_{k}(\psi^{\d^{*}(k)}_{q_{0},r})>$ as a function
$\o\rightarrow\C$.
Recall from Lemma 5 that
$$(14)\ \ \ \ \ <F^{\d^{*}(k)}_{k}(\psi^{\d^{*}(k)}_{q_{1},r})\vert L^{t}_{k\d^{*}(k)}\vert F^{\d^{*}(k)}_{k}(\psi^{\d^{*}(k)}_{q_{0},r})>=$$
$$(m/2\pi \hbar t)^{1/2}(2^{-2\d^{*}(k)}/\sqrt{N_{k}})$$
$$\int_{X_{k}^{**s}(y_{0},y_{1},r)}e^{i\Delta t(\sum_{j=1}^{n^{*}}S^{*}(x_{j-1},x_{j},\Delta t)/\hbar)}dx_{0}...x_{n^{*}},$$
where as before we have used the fact that $2^{-\d (k)}/\sqrt{N^{d}_{k}}=2^{-2\d (k)}/\sqrt{N_{k}}$.

Now recall that we write $P_{y_{0},y_{1},r}$ for the set of all paths $p\in P$
such that $y_{0}-(1/r)\le p^{k}_{0}<y_{0}+(1/r)$ and $y_{1}-(1/r)\le p^{k}_{n^{*}(k)}<y_{1}+(1/r)$
for all large enough $k$ ($r$ is a non-zero natural number and keep in mind that $y_{0}$, $y_{1}$ and $r$ determine what is large enough $k$
i.e. it does not depend on $p$).
We write $f=e^{iS(p)/\hbar}:P\rightarrow\C^{\o}$, for $S(p)$, see Section 2.
Now at least for all $k$  large enough,
if $X$ is a finite partition of $P_{y_{0},y_{1},r}$ of open sets with small enough mesh (compared to $k$ i.e. $<1/k$), by (14) and the definition of our integrals,
$$(15)\ \ \ \ \ <F^{\d^{*}(k)}_{k}(\psi^{\d^{*}(k)}_{y_{1},r})\vert L^{t}_{k\d^{*}(k)}\vert F^{\d^{*}(k)}_{k}(\psi^{\d^{*}(k)}_{y_{0},r})>=$$
$$(m/2\pi \hbar t)^{1/2}(2^{-2\d^{*}(k)}/\sqrt{N_{k}})(I_{f}(X)(k)).$$
Notice that this means that
$$(m/2\pi\hbar t)^{1/2}(2^{-2\d^{*}(k)}/\sqrt{N_{k}})\int_{P_{y_{0},y_{1},r}}e^{iS(p)/\hbar}dp$$
is the function
$$k\mapsto <F^{\d^{*}(k)}_{k}(\psi^{\d^{*}(k)}_{y_{1},r})\vert L^{t}_{k\d^{*}(k)}\vert F^{\d^{*}(k)}_{k}(\psi^{\d^{*}(k)}_{y_{0},r})>.$$

So by combining (13), (15) and Assumption 4 (as in the previous subsection) we get Theorem 1 i.e.
$$(16)\ \ \ \ \ K(y_{0},y_{1},t)=$$
$$(m/2\pi\hbar t)^{1/2}lim_{r\rightarrow\infty}lim_{D^{*}}(2^{-2\d^{*}(k)}r^{2}/4\sqrt{N_{k}})\int_{P_{y_{0},y_{1},r}}e^{iS(p)/\hbar}dp.$$
For $r^{2}/4$ in the formula, see $(9)$ in Subsection 3.1 or Lemma 5.

\medskip

{\bf Remark 5}. Let $P^{0}$ be the set of all continuously differentiable paths $p:[0,t]\rightarrow\R$.
Now the topology and the measure like function from $P$ induce a topology and
a measure like function to $P^{0}$ and these have a natural definition also in terms of $P^{0}$
alone, one does not need to define them via $P$. However there does not seem to be an obvious way to do the same for
our action like function $S$. But still one can, if one insists, write the formula (16) also using integration over $P^{0}$
(again one can compare this to Wiener measure).

\medskip

To get our formula (16) look more like the formula in the definition of Feynman path integral (the first version)
we would like to get rid of $r$ from the formula. However, here we have again the same problem that prevented us from
proving (12) in any other case than the free particle case and which forced us to make the Assumption 4,
again, see the discussion just before
the assumption.

\chapter{References}

[Fe] R. Feynman, The development of the space-time view of quantum electrodynamics, Nobel lecture, 1965,

\noindent
nobelprize.org/prizes/physics/1965/feynman/lecture.

[HH1] \AA. Hirvonen and T. Hyttinen, On eigenvectors, approximations and the Feynman propagator, Annals of Pure and Applied Logic, vol. 170, 2019, 109-135.

[HH2] \AA. Hirvonen and T. Hyttinen, On ultraproducts, the spectral theorem and rigged Hilbert spaces, Journal of Symbolic Logic, to appear.

[Me] J. Melo, Introduction to renormalization, University of Cambridge, 2019,

\noindent
arxiv.org/pdf/1909.11099.pdf

[Pe] D. Perepelitsa, Path integrals in quantum mechanics, manuscript, MIT, 2007,

\noindent
courses.physics.ucsd.edu/2018/spring/physics142/labs/lab6/pathintegrallecture.pdf.

[RS] M. Reed and B. Simon, Methods of Modern Mathematical Physics I: Functional analysis, Academic
Press, New York, 1972.

[Zi] B. Zilber, On model theory, non-commutative geometry and physics, draft, 2010.

\bigskip

Department of mathematics and statistics

University of Helsinki

Finland

\end